\documentclass[a4paper,11pt]{article}

\usepackage[hidelinks]{hyperref}
\usepackage{tikz,graphics}
\usepackage[latin1]{inputenc}
\usepackage{fancyhdr}
\usepackage{eurosym,amssymb,amsmath,amsthm,fullpage,graphicx}
\usepackage{lastpage}
\usepackage{setspace}
\usepackage{tabularx}
\usepackage{multicol,paralist}
\usepackage{fullpage,verbatim,enumerate}
\usepackage{mathptmx}

\usepackage[section]{placeins}

\usepackage{setspace}

\usepackage{multicol}
\usepackage{multirow}
\usepackage[compact]{title sec}
\usepackage{color}
\usepackage{colortbl}
\usepackage{xcolor}

\usepackage{hyperref}
\newtheorem{theorem}{Theorem}[section]
\newtheorem{corollary}[theorem]{Corollary}
\newtheorem{lemma}[theorem]{Lemma}
\newtheorem{prop}[theorem]{Proposition}
\newtheorem{fact}[theorem]{Fact}
\newtheorem{main}[theorem]{Main Result}

\theoremstyle{definition}\newtheorem{defi}[theorem]{Definition}
\newtheorem{con}[theorem]{Construction}
\newtheorem{remark}[theorem]{Remark}

\def\<{\langle}
\def\>{\rangle}

\newcommand{\cC}{\mathcal{C}}
\newcommand{\cM}{\mathcal{M}}

\newcommand{\diam}{\mathrm{Diam\,}}

\newcommand{\K}{\mathbb{K}}
\renewcommand{\L}{\mathbb{L}}

\newcommand{\cL}{\mathcal{L}}

\parskip\medskipamount

  \renewenvironment{thebibliography}[1]{%
    \begin{oldthebibliography}{#1}%
      \setlength{\parskip}{0ex}%
      \setlength{\itemsep}{0ex}%
  }%
  {%
    \end{oldthebibliography}%
  }

\title{On Exceptional Lie Geometries}
\author{\phantom{blablabla} Anneleen De Schepper \thanks{Department of Mathematics, Ghent University, Belgium, \texttt{Anneleen.DeSchepper@ugent.be}, supported by the Fund for Scientific Research Flanders|FWO Vlaanderen} \and Jeroen Schillewaert \thanks{Department of Mathematics, University of Auckland,New Zealand, \texttt{j.schillewaert@auckland.ac.nz}}\phantom{blablabla} \and Hendrik Van Maldeghem \thanks{Department of Mathematics, Ghent University, Belgium, \texttt{Hendrik.VanMaldeghem@ugent.be}} \and Magali Victoor \thanks{Department of Mathematics, Ghent University, Belgium, \texttt{Magali.Victoor@ugent.be}}}
 \date{Dedicated to the memory of Ernie Shult}

\begin{document}
\maketitle

\begin{abstract}
Parapolar spaces are point-line geometries introduced as a geometric approach to (exceptional) algebraic groups. We characterize a wide class of Lie geometries as parapolar spaces satisfying a simple intersection property. In particular many of the exceptional Lie incidence geometries occur. {In an appendix, we extend our result to the locally disconnected case and discuss the locally disconnected case of some other well known characterizations.} 
\end{abstract}


\section{Introduction}
\subsection{General context}
Buildings, sometimes more specifically called Tits-buildings, were introduced by Jacques Tits \cite{Tits:74} in order to have a geometric interpretation of semi-simple groups of algebraic origin (semi-simple algebraic groups, classical groups, groups of mixed type, (twisted) Chevalley groups). The definition of a building is, however, somewhat involved, and it does not immediately provide a good intuition. On the other hand, projective spaces and polar spaces|which are natural point-line geometries|can be given the structure of Tits-buildings and provide excellent permutation representations for the classical groups. In fact, projective and polar spaces are \emph{Grassmannians} of certain Tits-buildings. More exactly, since a Tits-building is a numbered simplicial complex, one can take all simplices of a certain type $T$ as point set, and then there is a well-defined mechanism that deduces a set of lines. The resulting point-line geometry is the so-called $T$-Grassmannian of the Tits-building.  For a certain choice of $T$, projective spaces and polar spaces emerge from spherical Tits-buildings of types $\mathsf{A}_n$ and $\mathsf{B}_n$, respectively. Other choices of $T$ for these and for other types of spherical Tits-buildings in general lead to \emph{parapolar spaces}, which were introduced and first studied by Cooperstein \cite{Co}. For a precise definition of these, see Section~\ref{Preliminaries}. 

Parapolar spaces are point-line geometries introduced to further simplify the geometric approach to Chevalley groups of exceptional types. The pure definition of a parapolar space is yet general enough to also capture many other geometries, in particular many Grassmannian geometries related to classical groups and (spherical) Tits-buildings. In general, the Grassmannian geometries related to a spherical Tits-building (or, as soon as the rank is at least 3, related to a semi-simple algebraic, classical or mixed group), and by extension also to a non-spherical one, are called \emph{Lie (incidence) geometries}.    A lot of work in the past went into characterizing, using additional properties, certain classes of parapolar spaces, preferably containing as many of exceptional type as possible. Among the main assumptions of these characterizations, one often finds the existence of  
a gap in the spectrum of dimensions of singular subspaces of symplecta arising from intersecting the latter with the perp of a point. 

{\subsection{Informal description of our Main Results}}

In the present paper we start with the observation that in the most popular exceptional Lie geometries, gaps also appear in another type of spectrum, namely, in the spectrum of the dimensions of the singular subspaces that occur as intersections of two symplecta. A parapolar space with at least one gap in that spectrum will be called here \emph{lacunary}, or, more exactly, $k$-lacunary, if $k$-dimensional singular subspaces never appear as the intersection of two symplecta, and all symplecta really possess $k$-dimensional singular subspaces.  The exceptionality of this behavior is then proved by our {Main Result \ref{main1}}, which completely classifies all lacunary parapolar spaces of symplectic rank at least 3. Moreover, we also include rank 2 symplecta, {see Main Result \ref{main2}}, 
using a harmless additional condition, which we explain in more detail below (Subsection~\ref{rank2}).

Our method of proof requires some local recognition results, in the style of ``if every point-residual of a parapolar space $\Omega$ is isomorphic to a certain fixed parapolar space, then the isomorphism type of $\Omega$ is unique and known''. We devote a section to such results, which seem to be new. Some of them can be proved using known characterizations, for others we use Tits' local approach to buildings \cite{Tits:83}. These results are, in our opinion, interesting in their own right and ready-made to apply in other situations. 


\subsection{Connection with the Freudenthal-Tits Magic Square}

In the 1960s Hans Freudenthal \cite{Fr} and Jacques Tits \cite{Tits1,Tits2} provided a remarkable uniform construction of some semi-simple Lie algebras, among which many of exceptional type. The types of these Lie algebras can be arranged in a symmetric $4\times 4$ square. In his habilitation thesis \cite{Tits3}, Tits describes a class of geometries related to Lie groups whose types are arranged in the same $4\times 4$ square. In fact, the geometries described by Tits are Grassmannians of certain buildings, most of which are parapolar spaces \emph{avant la lettre}. With the modern notation, the table of types looks as follows. 

\begin{table}[htbp]

\centering

\small

\begin{tabular}{| c | c | c | c |}\hline
$\mathsf{A}_{1,1}$  &   $\mathsf{A}_{2,\{1,2\}}$   &  $\mathsf{C}_{3,2}$ &  $\mathsf{F}_{4,4}$  \\ \hline
$\mathsf{A}_{2,1}$  &  $\mathsf{A}_{2,1}\times \mathsf{A}_{2,1}$ & $\mathsf{A}_{5,2}$ &  $\mathsf{E}_{6,1}$\\ \hline
$\mathsf{C}_{3,3}$  & $\mathsf{A}_{5,3}$ & $\mathsf{D}_{6,6}$ & $\mathsf{E}_{7,7}$\\ \hline
$\mathsf{F}_{4,1}$  &   $\mathsf{E}_{6,2}$ &  $\mathsf{E}_{7,1}$ & $\mathsf{E}_{8,8}$\\ \hline
\end{tabular}
\end{table}

Our approach allows to elegantly characterize many geometries (and their subgeometries) related to this so-called \emph{Freudenthal-Tits Magic Square}  as the $k$-lacunary parapolar spaces with symplectic rank at least $k+3$ (strong if $k=-1$); all geometries of the South-East $3\times3$ subsquare are thus captured. 

This is in fact our main motivation: single out the properties of the parapolar spaces in the Magic Square. Moreover these results will be of use in the investigation of the projective varieties associated to the square \cite{KSV,SVM,SVMbis}, in particular for the study of the Lagrangian Grassmannians (third row) and the adjoint varieties (fourth row).

Zooming in on the second row of the Magic Square, one obvious property of these parapolar spaces  is that every pair of symplecta intersects non-trivially. This corresponds to our notion of $(-1)$-lacunarity. Such (strong) parapolar spaces are classified in \cite{MinusOnePaper}.
That result, together with the classification of $0$-lacunary strong parapolar spaces admitting symplecta of rank 2,  is then taken as the first step of an inductive process to determine all lacunary parapolar spaces, see Subsection~\ref{structuur} below. 

\subsection{When some symplecta have rank 2} \label{rank2}
When there are symplecta of rank 2 around, then we 
assume that the parapolar space is strong
, i.e., every pair of intersecting lines is contained in either a singular subspace, or in a symplecton. The reason is that we consider the case that all symplecta have rank at least 3 as the main case, but the rank 2 case is needed in the induction process, where it only turns up under the strongness assumption. Hence, in this view, it is a harmless restriction. Note, however, that some important parapolar spaces escape in this way (e.g., the point-plane geometry arising from a projective space of dimension 3), but also some dull examples are excluded this way (e.g., the direct product of an arbitrary number of (possibly wild) polar spaces of rank 2).

A comment on the case of $0$-lacunarity is in order here. This notion expresses that in a parapolar space, symplecta that meet in a point, automatically share a line. In a private conversation with the second author, Shult mentioned that he would have liked to classify these parapolar spaces. This is also indicated by Exercise 13.26 in his book \cite{Shu:12}, which deals exactly with the $0$-lacunary hypothesis. However, what Shult was missing was the structure of the point-residuals; hence he was missing the classification of $(-1)$-lacunary parapolar spaces, in particular the ones with symplecta of rank 2. About the latter, though, he writes in \cite{Shultadvgeom} \emph{``It is not easy to live in a world with no symplecton of rank at least three in sight.''} Nevertheless, in our context,  we manage to deal with $0$-lacunary strong parapolar spaces whose symplecta are all of rank~2.

We end this subsection with another comment on Shult's work. It is tempting to use Shult's so-called \emph{Haircut Theorem} in \cite{Shulthaircut}, since (1) the hypotheses are residual, and (2) it is one of the rare results that does not assume constant rank. However, there seem to be problems with the final version of that paper as it appeared in \emph{Innovations in Incidence Geometry}, and so we chose not to use it, although a version of the Haircut Theorem using constant symplectic rank could be saved. Paraphrasing Shult's remark above, the moral is that \emph{``it is not easy to live in a world without the comfort of constant symplectic rank.''} Yet, in the present paper, we also manage to deal with such a world. 


\subsection{Structure of the paper}\label{structuur}

The paper is structured as follows. In \textbf{Section~\ref{MainResults}} we present our main results, and some corollaries. In \textbf{Section~\ref{Preliminaries}}, we provide all necessary definitions, and we review some known results in the theory of parapolar spaces, which we will use in our proof. In \textbf{Section~\ref{geometries}} we describe the examples and tabulate the geometries in the conclusions of our main results. \textbf{Section~\ref{localrecognition}} gathers the necessary local recognition lemmas. 

From then on, we proceed inductively. \textbf{Section~\ref{sec0lac}} deals with the $0$-lacunary case, and partially builds on the main result of \cite{MinusOnePaper}, namely, when all symps have rank at least 3. If some symps have rank 2, then we show that either we have a so-called imbrex geometry and we can use \cite{imbrex}, or the diameter of the parapolar space is $3$. In the latter case we prove the assumptions of a characterization result due to Kasikova and Shult \cite{KasikovaShult}. In \textbf{Section~\ref{sec1lac}}, we treat the $1$-lacunary case. The theoretic presence of symps of rank 2 is the only obstacle for directly combining our local recognition lemmas with the fact that the point residuals are $0$-lacunary parapolar spaces and hence classified in the previous section. However, we show that rank 2 symps cannot occur. From then on, the symplectic rank of the parapolar space is at least 3, and we can use an inductive approach, neatly supported by our local recognition lemmas. This is explained in detail in \textbf{Section~\ref{sec>2}}, concluding the proof of our main results.  

{
Finally, we refer to  \textbf{Appendix~\ref{ap}} for an extension of the current characterisation to the locally disconnected case. We present a formal and explicit method to do so, that can also be applied for other results concerning (locally connected) parapolar spaces, as we demonstrate by discussing some special cases (which also show the need of  formalising this procedure {beyond the principles noted down by Shult in~\cite{Shu:12}}).} 

\section{Main Results}\label{MainResults}
In this section we collect our main results. For the precise definitions and the notations we refer to Section~\ref{Preliminaries}; for the definitions and descriptions of the geometries in the conclusions, we refer to Section~\ref{geometries}. We content ourselves to informally define a \emph{$k$-lacunary} parapolar space as a parapolar with symplectic rank at least $k+1$ (i.e., all symplecta have rank at least $k+1$ as polar spaces) such that symplecta never intersect at precisely a $k$-dimensional singular subspace. If $k$ is not specified, we just say that the parapolar space is \emph{lacunary}. We say that a parapolar space has minimum symplectic rank $d$ if it has symplectic rank at least $d$ and there exists a symplecton of rank $d$. We also note that we do not consider polar spaces of infinite rank; these could be included by the interested reader with any definition. There does not seem to be a general agreement on what such polar spaces should precisely be, so we leave them out and consider inclusion of them as a minor additional effort that can be easily performed (they will arise in $(\mathsf{B})$ and $(\mathsf{D})$ below, but the rest remains unchanged).

Our main results are as follows. We combine them with the results from \cite{MinusOnePaper} in order to also include the $(-1)$-lacunary case.  For a more detailed version, mentioning the index $k$ for which the parapolar space is $k$-lacunary, we refer to Tables~\ref{lacmin1}~and~\ref{lac0}.

\begin{main}\label{main1} Let $\Omega=(X,\cL)$ be a locally connected  lacunary parapolar space with  symplectic rank at least $3$.  Then $\Omega$ is one of the following Lie incidence geometries ($\mathbb{K}$ is any commutative field, $\mathbb{L}$ is an arbitrary skew field):

\begin{compactenum}
\item[$\mathsf{(A)}$] $\mathsf{A_{5,3}}(\mathbb{L})$ or the line Grassmannian of a not necessarily finite-dimensional projective space of dimension at least $4$;

\item[$\mathsf{(B)}$] The line Grassmannian of an arbitrary thick polar space of rank at least $4$, that is, $\mathsf{B}_{n,2}(*)$;
\item[$\mathsf{(D)}$] $\mathsf{D}_{n,2}(\mathbb{K})$, $n\geq 4$, and a homomorphic image of $\mathsf{D}_{n,n}(\mathbb{K})$, $n\geq 5$ (isomorphic image if $n\leq 9$);
\item[$\mathsf{(E)}$] $\mathsf{E_{6,1}}(\mathbb{K})$, $\mathsf{E_{6,2}}(\mathbb{K})$, $\mathsf{E_{7,1}}(\mathbb{K})$, $\mathsf{E_{7,7}}(\mathbb{K})$, $\mathsf{E_{8,8}}(\mathbb{K})$, a homomorphic image of $\mathsf{E_{8,1}}(\mathbb{K})$, a homomorphic image of $\mathsf{E}_{n,1}(*)$, with $n\geq 9$ and $\mathsf{E}_n(*)$ any building of type $\mathsf{E}_n$;
\item[$\mathsf{(F)}$] $\mathsf{F_{4,1}}(*)$.
\end{compactenum}
\end{main}

If there are symps of rank 2, then local connectivity is not a suitable condition and we replace it by strongness, leading to the following theorem.

\begin{main}\label{main2}  Let $\Omega=(X,\cL)$ be a lacunary strong parapolar space with minimum symplectic rank $2$. Let $\Xi$ be the family of symps and let $\Sigma$ be the family of maximal singular subspaces (not assumed to be projective or finite-dimensional). 

Then $\Omega$ is one of the following geometries: 
\begin{compactenum}
\item[$\mathsf{(S)}$] \textbf{\em (Segre)} A Segre geometry $\mathsf{A_{2,1}}(\mathbb{*})\times \mathsf{A_{2,1}}(\mathbb{*})$;
\item[$\mathsf{(GQ)}$] \textbf{\em (Generalized Quadrangle)} An imbrex geometry such that $(X,\Sigma)$ is a generalized quadrangle and each member of $\Xi$ is an ideal subquadrangle of $(X,\Sigma)$;
\item[$\mathsf{(CP)}$] \textbf{\em (Cartesian Product)} the Cartesian product of a thick line with either an arbitrary linear space with thick lines, or an arbitrary polar space;
\item[$\mathsf{(DPS})$] \textbf{\em (Dual Polar Space)} a dual polar space of rank $3$, that is, $\mathsf{B_{3,3}}(*)$. 
\end{compactenum}
\end{main}

As a consequence we have:
\begin{corollary} \label{corMR}Let $\Omega$ be a locally connected parapolar space of minimum symplectic rank $d$.
 \begin{compactenum} 
\item[$(i)$] \label{xiv} If $d\geq 3$, then $\Omega$  is not $(d-1)$-lacunary.
\item[$(ii)$] \label{xv} If $d\geq 8$, then $\Omega$ is not lacunary.\end{compactenum} 
\end{corollary}\section{Preliminaries}\label{Preliminaries}
We refer to \cite{BuCo} and \cite{Shu:12} for the general background on parapolar spaces. To keep the paper self-contained and to fix notation, we provide a gentle introduction nonetheless. 

\subsection{Point-line geometries}\label{subsecdef}
\begin{defi} \em A pair $\Omega=(X,\cL)$ is a \emph{point-line geometry} if $X$ is a set and $\cL$ is a set of subsets of $X$ of size at least 2 covering $X$; the elements of $X$ are called {points} and those of $\cL$ lines. \end{defi}

Let $\Omega=(X,\cL)$ be a point-line geometry. Two distinct points $x,y$ of $X$ that are contained in a common line are called \emph{collinear}, denoted $x\perp y$. The set of points equal or collinear to a given point $x$ is denoted $x^\perp$, and for a set $S\subseteq X$, we denote $S^\perp=\bigcap_{s\in S}s^\perp$.  

A subset $Y\subseteq X$ is called a \emph{subspace of $\Omega$} if for every pair of collinear points $x,y\in Y$, all lines joining $x$ and $y$ are entirely contained in $Y$; it is called \emph{proper} if $Y\neq X$. A \emph{geometric hyperplane} is a proper subspace which intersects every line nontrivially. A subspace $Y\subseteq X$ is called \emph{singular} if every pair of distinct points of $Y$ is collinear.  For a subset $A$ of $X$, we denote by $\<A\>$ the intersection of all subspaces of $\Omega$ containing $A$. Then $\<A\>$ is a subspace itself, \emph{generated} by $A$; if the points of $A$ are pairwise collinear, then $\<A\>$ is a singular subspace. 

The \emph{collinearity graph} $\Gamma(X,\cL)$ of $\Omega$ is the graph on $X$ with collinearity as adjacency. A subspace $Y$ of $\Omega$ is called \emph{convex} if for every pair of points $x,y\in Y$, all points on any shortest path from $x$ to $y$ (in the collinearity graph) belong to $Y$. The intersection of all convex subspaces of $\Omega$ containing a given subset $A\subseteq X$ is called the \emph{convex closure} of $A$. A point-line geometry is called \emph{connected} if its collinearity graph is connected. The \emph{diameter} of a connected point-line geometry is the diameter of its collinearity graph.

\subsection{Parapolar spaces}

Before introducing parapolar spaces, we recall the definition of polar spaces (for more information, see Section 7.4 of \cite{BuCo}).

\begin{defi}\label{defpolarspaces}
A point-line geometry $\Delta=(X,\cL)$ is called a \emph{polar space} if the following axioms holds.
\begin{compactenum}
\item[(PS1)] Every line is \emph{thick}, that is, contains at least three points.
\item[(PS2)] No point is collinear to all other points.
\item[(PS3)] Every nested sequence of singular subspaces is finite.
\item[(PS4)] For any point $x$ and any line $L$, either one or all points on $L$ are collinear to $x$. 
\end{compactenum}
\end{defi}

Let $\Delta=(X,\cL)$ be a polar space. Then each of its singular subspaces is a projective space---that is, either arising from a vector space of dimension at least 4, an axiomatic projective plane, a thick line,  just one point, or the empty set---and hence its dimension is well-defined. Moreover, this also implies that each pair of collinear points define a unique line. One can show \cite{BuCo} that there is an integer $r \geq 2$ such that each singular subspace is contained in a singular subspace of dimension $r-1$, which we hence also refer to as a \emph{maximal singular subspaces}. We call $r$ the \emph{rank} of $\Delta$, denoted $\mathsf{rk}(\Delta)$. Note that (PS3) implies that the rank is finite, whereas some people also consider definitions of polar spaces allowing infinite rank.  We will however only consider polar spaces of finite rank, see Remark~\ref{finiterank}.


Crucial for the definition of parapolar spaces is the fact that $\Delta$ is the convex closure of any of its pairs of non-collineair points. Moreover, this convex closure can be obtained using a finite procedure, using the following sets: For non-collinear points $p,q$ of $\Delta$, we denote by $S_{p,q}$  the set $p^\perp \cap q^\perp$; for distinct collinear points $p,q $ of $\Delta$, we use the notation $S_{p,q}$ for the points on the line $pq$.

\begin{fact}\label{convexclosure} Let $\Delta=(X,\cL)$ be a polar space and let $p,q \in X$ be non-collinear. Then the convex closure of $p$ and $q$ consists of all points in $\{S_{p',q'} \mid p' \in \<p,r\>, q' \in \<q,r\>, r \in S_{p,q}\}$, and coincides with $\Delta$. \end{fact}

\begin{defi}\label{defparapolarspaces} A point-line geometry $\Omega=(X,\cL)$ is called a {\em parapolar space} if the following hold:
\begin{compactenum}
\item[(PPS1)] $\Omega$ is connected and, for each line $L$ and each point $p \notin L$, $p$ is collinear to either none, one or all of the points of $L$ and there exists a pair $(p,L)\in X \times \cL$ such that $p$ is collinear to no point of $L$. 
\item[(PPS2)] For each pair of non-collinear points $p$ and $q$ in $X$ with $p^\perp \cap q^\perp$ non-empty, either
\begin{compactenum}[(a)]
\item the convex closure of $\{p,q\}$ is a polar space, called a \emph{symplecton} or briefly a \emph{symp}; or
\item $p^\perp \cap q^\perp$ is a single point, and then $p$ and $q$ are called \emph{special}.
\end{compactenum}
\item[(PPS3)] Every line is contained in at least one symplecton.
\end{compactenum}
\end{defi}

Let $\Omega=(X,\cL)$ be a parapolar space.  We call $\Omega$ \emph{strong} if there are no special pairs; and we denote the set of symps by $\Xi$. Let $d(\Xi):=\{\mathsf{rk}(\xi) \mid \xi \in \Xi\}$.   We say that $\Omega$ has \emph{symplectic rank at least $d$} if $\min d(\Xi) \geq d$;  \emph{minimum symplectic rank $d$} if $\min d(\Xi) =d$ and  \emph{(uniform) symplectic rank $d$} if $d(\Xi)=\{d\}$. 



\begin{fact}\label{planelinesymp} If all points of a line $L$ contained in a symp $\xi$ of rank at least $3$ are collinear to a point $p$, then $p$ and $L$ are contained in a symp and hence generate a projective plane. Consequently, if the symplectic rank is at least $3$, each singular subspace is a projective space.\end{fact}


\begin{fact}\label{subspacesymp}
Let $\Omega$ be a parapolar space of uniform symplectic rank $d$. Then every singular subspace of dimension at most $d-1$ is contained in some symp.
\end{fact}

\subsection{The local structure of a parapolar space}





\begin{defi}\label{lc} Let  $\Omega=(X,\cL)$ be a parapolar space and $p$ one of its points. We call $\Omega$ \emph{locally connected at $p$} if each two lines through $p$ are contained in a finite sequence of singular planes consecutively intersecting in lines through $p$; $\Omega$  is called locally connected if it is locally connected at $p$ for all $p\in X$.
\end{defi}

\begin{defi}\label{point-residue} Let $\Omega=(X,\cL)$ be a parapolar space and let $p$ be one of its points. We define the \emph{point-residual at $p$}, denoted $\Omega_p=(X_p,\cL_p)$, as follows:
\begin{compactenum}[$-$]
\item $X_p$ is the set of lines through $p$;
\item $\cL_p$ is the set of planar line-pencils with vertex $p$ contained in singular planes through $p$ which are contained in a symp of $\Omega$.
\end{compactenum}
\end{defi}
The reason that we only consider singular planes which are contained in symps is that we want Axiom~(PPS3) to hold in $\Omega_p$; note that this is only a minor requirement since Fact~\ref{planelinesymp} implies that each singular plane is contained in a symp if $\Omega$ has symplectic rank is at least 3. Suppose $\xi\in \Xi$ contains $p$. If $\mathsf{rk}(\xi) \geq 3$, then $\xi$ corresponds to a polar space $\xi_p$ in $\Omega_p$ of rank $\mathsf{rk}(\xi)-1$, whose points are the lines of $\xi$ containing $p$ and whose lines are the planar line-pencils with vertex $p$ in $\xi$; If $\mathsf{rk}(\xi)=2$, then $\xi$ ``vanishes'' in $\Omega_p$.

The following fact is based on Theorem 13.4.1 of \cite{Shu:12}.

\begin{fact}\label{pt-pps} Let $\Omega=(X,\cL)$ be a parapolar space, assumed to be strong if its minimum symplectic rank is $2$, and let $p$ be any of its points.  Let $C$ be a connected component of $\Omega_p$. Then either:
\begin{compactenum}[$-$]
\item $C$ is a single element of $\cL$ (which then corresponds to a line of $\Omega$ through $p$ only contained in symps of rank $2$);
\item $C$ corresponds to the lines through a point in a symp $\xi$ of $\Omega$ of rank at least $3$ (which happens if no line of $\xi$ is contained in a second symp of $\Omega$ of rank at least $3$);
\item $C$ is a strong parapolar space. There is a bijective correspondence between the singular subspaces of $\Omega$ through $p$ and the singular subspaces of $\Omega_p$ and between the symps of $\Omega$ through $p$ of rank at least $3$ and the symps of $\Omega_p$.
\end{compactenum}
\end{fact}

The following fact is almost by definition:

\begin{fact}\label{conn}
Let $\Omega=(X,\cL)$ be a parapolar space with symplectic rank at least $3$.
For all $p\in X$, $\Omega$ is locally connected at $p\in X$ if and only if $\Omega_p$ is connected.
\end{fact}




\subsection{Useful theorems}

We end this section by listing two theorems which we will use during our classification. Firstly, Kasikova and Shult obtained the following characterization result \cite[Theorem 2]{KasikovaShult}. Note that they labelled parapolar spaces differently,
we have phrased their theorem here conforming to our convention of Bourbaki labelling, and using the notation we introduced above (in particular we only consider parapolar spaces with symps of finite rank).

\begin{theorem}\label{KasikovaShult2}
Suppose $\Gamma=(X,\cL)$ is a strong parapolar space with these three properties:
\begin{compactenum}
\item For every $x\in X$ and $\xi\in \Xi$, the set  $x^\perp \cap \xi$ is non-empty;
\item The set of points at distance at most 2 from any point $x$ forms a geometric hyperplane of $\Gamma$;
\item If the symplectic rank is at least 3,  every maximal singular subspace is finite-dimensional.
\end{compactenum}
Then $\Gamma$ is one of the following:
\begin{compactenum}[$\bullet$]
\item $\mathsf{D_{6,6}}(\K)$, $\mathsf{A_{5,3}}(\K)$, or $\mathsf{E_{7,1}}(\K)$ (where $\K$ is any field),
\item a dual polar space of rank $3$,
\item a product geometry $L\times \Delta$, where $L$ is a thick line and $\Delta$ is a polar space of rank at least 2.
\end{compactenum}
\end{theorem}

We also use some results about imbrex geometries, see \cite{imbrex}. An \emph{imbrex geometry} is a 
strong parapolar space of diameter $2$ and symplectic rank $2$ such that every pair of symps which share a point $p$ and both intersect a certain line at points not collinear to $p$, share a line.
\begin{prop}\label{prop:imbrex}
Let $\Omega=(X,\cL)$ be an imbrex geometry and let $\Sigma$ be the set of maximal singular subspaces of $\Omega$. Then either $\Omega$ is isomorphic to a product space $Y\times Z$, where $Y$ and $Z$ are arbitrary linear spaces with thick lines, or $(X,\Sigma)$ is a thick generalized quadrangle and each element of $\Xi$ is an ideal subquadrangle of $(X,\Sigma)$; moreover no maximal singular subspace is a projective space.  
\end{prop}

\section{Description of the (Lie incidence) geometries}\label{geometries}

Most of the geometries in the conclusion of  Main Result~\ref{main1}  are Lie incidence geometries, i.e., they arise from Tits-buildings or their quotients, by selecting a node $k$ of the corresponding Coxeter diagram and considering the so-called \emph{$k$-Grassmannian} of the (quotient of the) building. Below we give an overview of these $k$-Grassmannians, after introducing some generalities and terminology in the theory of buildings. 

\subsection{Tits-buildings and their $k$-Grassmannians}
Most spherical buildings arise, as the geometry of the Borel subgroup and the parabolic subgroups, from simple algebraic groups and close relatives like classical groups, (twisted) Chevalley groups and groups of mixed type. 
More exactly, we can think of a building of rank $n$ as an $n$-partite graph $G$ where each partition class consists of the cosets of a parabolic subgroup (with respect to a fixed chosen Borel subgroup) and where adjacency is given by ``intersecting non-trivially''.  The $n$-cliques are called \emph{chambers} and each clique is contained in some chamber. When each clique of size $n-1$ is contained in at least three chambers, then the building is called \emph{thick}. 

The induced bipartite graph obtained from an $(n-2)$-clique $C$ by considering all vertices $v$ such that $C\cup\{v\}$ is an $(n-1)$-clique is a building of rank 2. These graphs have diameter $\ell$ and girth $2\ell$, for some natural number $\ell\geq 2$. If $\ell=2$, this is a complete bipartite graph; if $\ell=3$, then this is the incidence graph of a projective plane; if $\ell=4$, then we have the incidence graph of a generalized quadrangle (in case there are at least three points per line, this is a polar space of rank 2). We will not need $\ell\geq 5$ here.

 It so happens that $\ell$ only depends on the $n-2$ partition classes containing a vertex of the $(n-2)$-clique $C$. So we can build a diagram with set of nodes the partition classes  of $G$ and have no edge, a single edge or a double edge between two nodes if $\ell=2,3,4$, respectively, where $C$ is any $(n-2)$-clique not containing any member of the partition classes corresponding to the two nodes under consideration. This is the \emph{Coxeter diagram} of the building. There is an enhanced notion of Dynkin diagram, but we will not need this. We number the nodes of the Coxeter diagram according to the Bourbaki labeling \cite{Bourbaki}. 

\begin{defi}\label{kgras} For a certain node $k$, the $k$-Grassmannian geometry is then the point-line geometry with set of points the partition class of $G$ corresponding to the node $k$; and a  line is the set of vertices of type $k$ contained in a chamber with a given $(n-1)$-clique not containing a vertex of type $k$. \end{defi}

In general, for a building of type $\mathsf{X}_n$, its $k$-Grassmannian is a Lie (incidence) geometry, and we refer to it by its Coxeter type $\mathsf{X}_{n,k}$. 
Below we give an overview of the Lie incidence geometries that we will encounter in this paper. For details about the corresponding parapolar spaces, such as symplectic rank, singular rank, strongness and diameter, we refer to Subsection~\ref{tabellekes}.

\paragraph*{(A) Grassmannians of vector spaces|Projective spaces}

Let $\Delta$ be a thick building of type $\mathsf{A}_n$ for $n\in \mathbb{N}$, i.e.,  $\Delta$ corresponds to a projective space $\mathbb{P}$ of dimension $n$. For any $k \in\mathbb{N}$ with $1 \leq k \leq n$, we consider the $k$-Grassmannian of $\Delta$. By Definition~\ref{kgras}, its points are the $(k-1)$-dimensional subspaces of $\mathbb{P}$, and a line is the set of such $(k-1)$-spaces  containing a fixed $(k-2)$-space and  contained in a fixed $k$-space.   If $\mathbb{P}$ is coordinatised by a skew field $\mathbb{L}$, then we denote this Grassmannian by $\mathsf{A}_{n,k}(\mathbb{L})$; if this is not necessarily  the case (i.e., if $n\leq 2$), then by $\mathsf{A}_{n,k}(*)$.

If $k\in\{1,n\}$, we just obtain $\mathbb{P}$ or its dual. If $n\geq 4$ and $k\notin\{1,n\}$, $\mathsf{A}_{n,k}(\L)$ is a strong parapolar space. The fact that the dual of $\mathbb{P}$ is a projective space of dimension $n$, allows us to restrict our attention to values of $k$ in $\{2,..., \frac{n}{2}\}$.  We add that these Grassmannians can be defined completely similarly for $n=\infty$, where $\infty$ denotes any cardinal.

\paragraph*{(B) Thick polar spaces and polar Grassmannians}

Let $\Delta$ be a thick building of type $\mathsf{B}_n$ for $n \in \mathbb{N}$ with $n\geq 2$, i.e.,  $\Delta$ corresponds to a \emph{thick} polar space $\Gamma$ of rank $n$; thick meaning that for each singular subspace of dimension $n-2$, there are at least three maximal singular subspaces containing it. For any $k \in \mathbb{N}$ with $1\leq k\leq n$, we  consider the $k$-Grassmannian of $\Delta$. By Definition~\ref{kgras}, its points are the $(k-1)$-dimensional singular subspaces of $\Delta$, and a line is a set of such $(k-1)$-spaces containing a fixed $(k-2)$-space and, if $k<n$, contained in a fixed singular subspace of dimension $k$.  Since a polar space is usually not uniquely defined by the underlying skew field of the singular projective spaces, we denote its $k$-Grassmannian by $\mathsf{B_{n,k}}(*)$. 

If $k=1$, we  obtain $\Gamma$; moreover,  $\mathsf{B_{2,1}}(*)$ and $\mathsf{B_{2,2}}(*)$ are  generalised quadrangles (they are each other's dual). If $n >2$ and $k >1$, then $\mathsf{B}_{n,k}(*)$ is  a parapolar space;  if $k=n$ it is called a \emph{dual polar space}.


\paragraph*{(D) Non-thick polar spaces and their Grassmannians}
Let $\Delta$ be a thick building of type $\mathsf{D}_n$ for $n \in \mathbb{N}$ with $n\geq 3$, i.e.,  $\Delta$ corresponds to a \emph{non-thick} polar space $\Gamma$ of rank $n$; non-thick meaning that for each singular subspace of dimension $n-2$, there are exactly two maximal singular subspaces containing it. Then there are two natural families $\mathcal{M}_1$ and $\mathcal{M}_2$ of maximal singular subspaces in $\Gamma$. Moreover, $\Gamma$ is defined over a (commutative) field $\K$. For any $k \in \mathbb{N}$ with $1\leq k\leq n$, we  consider the $k$-Grassmannian of $\Delta$, which we denote by $\mathsf{D}_{n,k}(\mathbb{K})$. 

When $k \leq n-2$, this is completely similar to the above case where $\Gamma$ is a thick polar space, so suppose $k=n-1$ or $k=n$. Again according to Definition~\ref{kgras}, the point set of the $k$-Grassmannian then coincides with $\mathcal{M}_1$ or $\mathcal{M}_2$; and a line is the set of such subspaces containing a fixed subspace of dimension $n-3$. The symmetry of the diagram implies that $\mathsf{D}_{n,n-1}(\mathbb{K})$ is isomorphic to $\mathsf{D}_{n,n}(\mathbb{K})$, and we will always use the latter notation. If $n=3$, we obtain a projective space of dimension 3; if $n=4$ then we obtain a polar space isomorphic to $\Gamma$; and if $n \geq 5$ we obtain a strong parapolar space, also called a \emph{half spin geometry}. 
Also here, the case $k=2$ corresponds to the long root geometries.

\paragraph*{(E) Exceptional parapolar spaces of type $\mathsf{E}_i$, $i=6,7,8$}
Let $\Delta$ be a thick building of type $\mathsf{E}_n$ for $n\in\{6,7,8\}$. Then $\Delta$ is defined over the field $\mathbb{K}$. For any $k\in\mathbb{N}$ with $ 1 \leq k \leq n$, the $k$-Grassmannian (cf.\ Definition~\ref{kgras}) is denoted by $\mathsf{E}_{n,k}(\K)$. There  are certain choices for $k$ such that  $\mathsf{E}_{n,k}(\K)$ has small diameter and constant symplectic rank. We record that $\mathsf{E_{6,1}}(\mathbb{K})$ has diameter 2, $\mathsf{E_{7,7}}(\mathbb{K})$, $\mathsf{E_{7,1}}(\mathbb{K})$, $\mathsf{E_{8,8}}(\mathbb{K})$ have diameter 3 and $\mathsf{E_{8,1}}(\mathbb{K})$ has diameter 5.

\paragraph*{(F) Metasymplectic spaces} 
Let $\Delta$ be a building of type $\mathsf{F}_4$, i.e.,  $\Delta$ corresponds to a \emph{metasymplectic space}, after Freudenthal \cite{Fr2}. This metasymplectic space is the $1$-Grassmannian of $\Delta$ and is a non-strong parapolar space of diameter 3.  Also, $\Delta$ is thick if and only if the corresponding metasymplectic space has no non-thick symps. A thick building of type $\mathsf{F}_4$ is  determined by a certain pair of division rings, but there is no need to explain this in more detail, so we just denote the corresponding metasymplectic spaces by $\mathsf{F_{4,1}}(*)$.  By the symmetry of the diagram, this class coincides with the $4$-Grassmannians $\mathsf{F_{4,4}}(*)$ of $\Delta$ (which does not mean that the $1$-Grassmannian of $\Delta$ is isomorphic to its $4$-Grassmannian though). We will not need other $k$-Grassmannians of $\Delta$. 

\paragraph*{Parapolar spaces from non-spherical buildings} Let $\Delta$ be a building corresponding to the diagram $\mathsf{E_n}$ with $n\geq 9$, obtained from $\mathsf{E_8}$ by extending it at the long arm with a simple path of length $n-8$, extending the Bourbaki numbering in the obvious way.  The case $n=9$ is exactly $\widetilde{\mathsf{E}}_8$. These buildings need not be unique given a field, so we denote their $1$-Grassmannians by $\mathsf{E}_{n,1}(*)$. These are non-strong parapolar spaces with unbounded diameter.

\paragraph*{Homomorphic images of the $k$-Grassmannians}
If the diameter of the parapolar space in one of the examples above is at least 5, then there might exist homomorphic images that are again parapolar spaces. Indeed, consider for example the $n$-Grassmannian $\mathsf{A}_{2n-1,n}(\mathbb{L})$, $n\geq 5$, arising from a vector space $V$, and suppose that $\sigma$ is a polarity of $V$ of Witt index at most $n-5$. Identifying  the $n$-spaces of $V$ that correspond to each other under $\sigma$, produces a homomorphic image  that is again a (locally connected) parapolar space with the same local properties as $\mathsf{A}_{2n-1,n}(\mathbb{L})$. We write a superscript $h$ to indicate a homomorphic (possibly isomorphic) image. In the above examples that are relevant to us, homomorphic images are possible in the geometries $\mathsf{A}_{2n-1,n}(\mathbb{L})$, $n\geq 5$, $\mathsf{D}_{n,n}(\mathbb{K})$, $n\geq 9$, and $\mathsf{E}_{n,1}(\mathbb{*})$, $n\geq 8$. 

\paragraph*{Long root geometries among the $k$-Grassmannians}
From the above Grassmannians,  $\mathsf{B}_{n,2}(*)$ (for $n \geq 3$), $\mathsf{D}_{n,2}(*)$ (for $n \geq 4$), $\mathsf{E}_{6,2}(\K)$, $\mathsf{E}_{7,1}(\K)$, $\mathsf{E}_{8,8}(\K)$ and $\mathsf{F}_{4,1}(*)$  are instances of so-called \emph{long root geometries}. The exact meaning of that is not essential for this paper, but we do mention that it implies that  these parapolar space have diameter $3$ and are non-strong, and hence cannot appear as point-residuals in another parapolar space (and that is most relevant for the present paper).

\paragraph*{Cartesian product spaces} 
Let $\Omega_i=(X_i,\cL_i)$, $i=1,2$, be two point-line geometries. Define the Cartesian product space $\Omega:=\Omega_1\times \Omega_2$ as the point-line geometry with point set the Cartesian product $X_1\times X_2$ and whose lines are  $\{p_1\} \times L_2$, for $p_1 \in X_1$ and $L_2\in \mathcal{L}_2$ and $L_1 \times \{p_2\}$, for $L_1 \in \mathcal{L}_1$ and $p_2 \in X_2$. 

Parapolar spaces that we encounter in our Main Result arising as a Cartesian product are those where $\Omega_1$ is a thick line (i.e., a set of at least three points, also referred to as $\mathsf{A}_{1,1}(*)$) and $\Omega_2$ is either a projective space (or more generally a \emph{linear space}, i.e., a point-line geometry in which each two points are on a unique line, with thick lines) or a polar space of rank at least 2, or where both  $\Omega_1$ and $\Omega_2$ are projective spaces, say of respective dimensions $n,m$ (with $nm>1$). In the latter case, $\Omega$ is the \emph{Segre geometry of type $(n,m$)}.

In general, if $\Omega_1$ and $\Omega_2$ are parapolar spaces, with set of symps $\Xi_1$ and $\Xi_2$ respectively, then $\Omega$  is again a parapolar space, the symps of which are formed by $\{p_1\} \times \xi_2$ for $p_1\in X_1$ and $\xi_2\in \Xi_2$, likewise $\xi_1 \times \{p_2\}$ for $\xi_1\in \Xi_1$ and $p_2\in X_2$, but also $L_1 \times L_2$ for $L_1\in \mathcal{L}_1$ and $L_2\in\mathcal{L}_2$.
If both $\Omega_1$ and $\Omega_2$ are strong, then also $\Omega$ is strong. Its diameter equals $\diam\Omega_1+\diam\Omega_2$, and it always has minimum symplectic rank~$2$. 




\subsection{The Main Result in tabular form}\label{tabellekes}

We now tabulate all $k$-lacunary parapolar spaces $\Omega$ which are either of minimum symplectic rank 2 (and hence strong by assumption) or which are locally connected. We include the following parameters: the rank and thickness of the symplecta, structure and dimensions of the maximal singular subspaces, diameter and strongness of $\Omega$.  
As the tables below will show, these parameters suffice to distinguish between all $k$-lacunary parapolar spaces with projective singular subspaces, for a given $k$. This is crucial to our inductive approach, as we need to know that for any two points $x,x' \in \Omega$, the point-residuals $\Omega_x$ and $\Omega_{x'}$ are isomorphic. 
As already hinted at in Subsection~\ref{structuur}, there are two natural classes to consider: those with symplectic rank $d \geq k+3$ (because they stem from $(-1)$-lacunary parapolar spaces) and $d=k+2$ (because they stem from $0$-lacunary parapolar spaces of minimum symplectic rank 2). 

\paragraph*{The $k$-lacunary parapolar spaces with (minimum) symplectic rank $d \geq k+3$.}
The $k$-lacunary parapolar spaces $\Omega$ in Table~\ref{lacmin1}  have minimum symplectic rank $d \geq k+3$. The table is such that, for any cell in any row, going one cell to the left in that row yields its point-residual. The set $S$ displays the dimensions of the maximal singular subspaces. White cells contain strong parapolar spaces of diameter 2, grey cells contain the ones of diameter 3. If the parapolar space appears in white, it means that it is non-strong. Since taking a point-residual yields a strong parapolar space, only the rightmost parapolar spaces can be non-strong. For each of the parapolar spaces in Table~\ref{lacmin1}, the symplecta are non-thick.

\begin{table}[h]
\begin{center}
\resizebox{\textwidth}{!}{
\begin{tabular}{|c||c||c|c|c|c|c|c|}
\hline

$d$ 					&$S$&$k=-1$								& $k=0$								& $k=1$						& $k=2$						& $k=3$ 						& $k=4$ 						  \\ \hline\hline

\multirow{2}{*}{$k+3$}	&$\{k+2,k+3\}$ &$\mathsf{A}_{1,1}(*) \times \mathsf{A}_{2,1}(*)$ 	& $\mathsf{A}_{4,2}(\L)$					& $\mathsf{D}_{5,5}(\K)$			& $\mathsf{E}_{6,1}(\K)$			& \cellcolor[gray]{0.8} $\mathsf{E}_{7,7}(\K)$	& \cellcolor[gray]{0.8} $\textcolor{white}{\mathsf{E}_{8,8}(\K)}$          	  \\ \cline{2-8}
					&$\{k+3\}$ &$\mathsf{A}_{2,1}(*) \times \mathsf{A}_{2,1}(*)$  	& \cellcolor[gray]{0.8}$\mathsf{A}_{5,3}(\L)$	& \cellcolor[gray]{0.8} $\textcolor{white}{\mathsf{E}_{6,2}(\K)}$	&					&							&	 					  \\ \hline\hline

\multirow{2}{*}{$k+4$}	&$\{k+3,k+4\}$ &$\mathsf{A}_{4,2}(\L)$ 					& $\mathsf{D}_{5,5}(\K)$					& $\mathsf{E}_{6,1}(\K)$			&  \cellcolor[gray]{0.8} $\mathsf{E}_{7,7}(\K)$	&\cellcolor[gray]{0.8} $\textcolor{white}{\mathsf{E}_{8,8}(\K)}$           		&	     					  \\ \cline{2-8}
					&$\{k+3,k+5\}$ &$\mathsf{A}_{5,2}(\L)$						& \cellcolor[gray]{0.8}$\mathsf{D}_{6,6}(\K)$	&  \cellcolor[gray]{0.8} $\textcolor{white}{\mathsf{E}_{7,1}(\K)}$	&  							&					       		&	     				          \\ \hline\hline
					
$k+6$				&$\{k+5,k+6\}$ &$\mathsf{E}_{6,1}(\K)$						&  \cellcolor[gray]{0.8}$\mathsf{E}_{7,7}(\K)$	& \cellcolor[gray]{0.8} $\textcolor{white}{\mathsf{E}_{8,8}(\K)}$	&  							&					       		&	     					  \\ \hline
\end{tabular}}
\end{center}
\vspace{1em}
\caption{The $k$-lacunary parapolar spaces with symplectic rank $d \geq k+3$.\label{lacmin1}}
\end{table}
\paragraph*{The $k$-lacunary parapolar spaces with minimum symplectic rank $k+2$.}
The $k$-lacunary parapolar spaces $\Omega$ in Table~\ref{lac0}  have minimum symplectic rank $d=k+2$.  This time, the symps can be thick; but thick and non-thick symps only occur in the same parapolar spaces in  the second last row of Table~\ref{lac0}. Also, only in this row and the row before, the rank is not uniform. For both reasons, we call the symps ``mixed'' in the second last row. We will use some abbreviations:
\begin{compactenum}
\item[$(\mathsf{GQ})$] This refers to the case where $\Omega$ is an imbrex geometry and $(X,\Sigma)$ is a generalized quadrangle (hence GQ) and each member of $\Xi$ is an ideal subquadrangle of $(X,\Sigma)$.
\item[$(\mathsf{LS})$] Stands for a linear space with only thick lines and which contains at least two lines.
\item[($\cdot^h$)] Recall that, if $Z$ is some locally connected parapolar space, then $Z^h$ denotes a locally connected homomorphic image of it, which is still a parapolar space (fibers only contain points at distance at least 5 from each other).
\item[($\mathsf{A}_1$)] Shorthand for $\mathsf{A}_{1,1}(*)$.
\end{compactenum}
The various $\mathsf{A_1}$ (thick lines) appearing in the same cell need not be isomorphic, i.e., they can have different sizes. 

Again, $S$ is the set of dimensions of the maximal singular subspaces of $\Omega$. If a class of parapolar spaces is such that its maximal singular subspaces can \emph{never} be projective, then we say that this dimension is empty (``$\emptyset$'').  If a class of parapolar spaces is such that its maximal singular subspaces are \emph{in some, but not all, cases} projective, we write ``$-$'', because in case the maximal singular subspaces are not projective, their dimension is undefined (by definition, this notation implies that, in case the maximal singular subspace happens to be projective, its dimension is at least 2). The number $n$ denotes the rank of the corresponding building.  The shades of grey reflect the diameter: the darker the cell, the higher the diameter; the white color indicates that $\Omega$ is non-strong. Also here, reading from right to left corresponds to taking point-residuals; so non-strong parapolar spaces can only occur at the rightmost position in a row. \begin{table}[h]
\begin{center}
\resizebox{\textwidth}{!}{
\begin{tabular}{|c||c||c|c|c|c|c|c|}
\hline
symps			&$S$	 	 &$k=0$					& $k=1$						& $k=2$					& $k=3$					& $k=4$ 					& $k=5$ 					  \\ \hline\hline
\multirow{6}{*}{non-thick} &
\multirow{4}{*}{\parbox[c]{1,7cm}{$\{1,-\}$\phantom{cc} (if $k=0$)\\ \mbox{} \\ $\{k\!+\!\!1\!,n\!\!-\!\!1\!\}$ (if $k\geq 1$)}} 		&\multirow{4}{*}{$\mathsf{A}_{1} \times \mathsf{LS}$}	&\multirow{4}{*}{\parbox[c]{1cm}{$\mathsf{A}_{n,2}(\L)$\\ $\phantom{ci}_{n\geq 4}$}}\phantom{i}
																			&$\mathsf{D}_{5,5}(\K)$		& ${\mathsf{E}_{6,1}(\K)}$&  \cellcolor[gray]{0.8} $\mathsf{E}_{7,7}(\K)$& \cellcolor[gray]{0.8}  $\textcolor{white}{\mathsf{E}_{8,8}(\K)}$  \\ \cline{5-8}
				      &	& 								& 						&\cellcolor[gray]{0.8}  ${\mathsf{D}_{6,6}(\K)}$			&  \cellcolor[gray]{0.8} $\textcolor{white}{\mathsf{E}_{7,1}(\K)}$& 					 	 &  \\ \cline{5-8}
				      && 								&						&\cellcolor[gray]{0.8}${\mathsf{D}_{7,7}(\K)}$			&  \cellcolor[gray]{0.5} $\textcolor{white}{{\mathsf{E}_{8,1}(\K)^h}}$&  						 &   \\ \cline{5-8}
				      &&								&						&\cellcolor[gray]{0.65} \raisebox{0,2em}{\parbox[c]{1cm}{$\mathsf{D}_{n,n}(\K)^h$\\[-6pt] \phantom{cc}$_{n \geq 8}$}	\phantom{x.}}		& 	\cellcolor[gray]{0.5} \raisebox{0,2em}{\parbox[c]{1cm}{$\textcolor{white}{{\mathsf{E}_{n,1}(*)^h}}$\\[-6pt] \phantom{cc}$_{n \geq 9}$}	\phantom{x.}}									&  				 		 &   \\\cline{2-8}
&	$\emptyset$		&$\mathsf{GQ}$					&							&						&						&						&							
				     	  \\ \cline{2-8}
				    
& $\{k+1\}$		&   \cellcolor[gray]{0.8}$\mathsf{A}_{1} \times \mathsf{A}_{1} \times \mathsf{A}_{1}$&  \cellcolor[gray]{0.65}$\textcolor{white}{\mathsf{D}_{4,2}(\K)}$&			&						&						&	  

 \\ \cline{2-8}
				    
& $\{k+1, n+k-2\}$		&   \cellcolor[gray]{0.8}\raisebox{0,2em}{\parbox[c]{2,3cm}{$\mathsf{A}_{1}\! \!\times\! \mathsf{D}_{n-1,1}(\mathbb{K})$\\[-6pt] \phantom{cccc}$_{n \geq 4}$}}&  \cellcolor[gray]{0.65} \raisebox{0,2em}{\parbox[c]{1cm}{\textcolor{white}{$\mathsf{D}_{n,2}(\mathbb{K})$}\\[-6pt] \phantom{cc}$_{n \geq 5}$}}	\phantom{c}&			&						&						&	  

\\ \hline\hline			  

mixed	     & $\{k+1,n+k-2\}$	&  \cellcolor[gray]{0.8} \raisebox{0,2em}{\parbox[c]{2,3cm}{$\mathsf{A}_{1} \times \mathsf{B}_{n-1,1}(*)$\\[-6pt] \phantom{cccc}$_{n \geq 3}$}}	
				& \cellcolor[gray]{0.65} \raisebox{0,2em}{\parbox[c]{1cm}{\textcolor{white}{$\mathsf{B}_{n,2}(*)$}\\[-6pt] \phantom{cc}$_{n \geq 4}$}}	\phantom{c}&						&						& 						&  	  \\ \hline\hline
thick 	 & $\{k+1\}$ &  \cellcolor[gray]{0.8} $\mathsf{B}_{3,3}(*)$		& \cellcolor[gray]{0.65}  $\textcolor{white}{\mathsf{F}_{4,1}(*)}$			&						&						&						&  \\ \hline

\end{tabular}}
\end{center}
\vspace{1em}
\caption{The $k$-lacunary parapolar spaces with minimum symplectic rank $k+2$.\label{lac0}}
\end{table}

It is not a coincidence that there are never three grey cells in one row, and that a cell written with white letters is always preceded by a grey cell. This follows from the property that a $\Omega$ is strong and of symplectic rank at least 3 if and only if its point-residuals have diameter 2. 

\subsection{Non-examples|weakly lacunary parapolar spaces}

We motivate our requirement that $k$-lacunary parapolar spaces only contain symps of rank at least $k+1$. It is obvious that at least two symps have rank $k+1$, since otherwise the $k$-lacunarity is no restriction at all. So let us call a parapolar space \emph{weakly $k$-lacunary} if no pair of symps intersects in a $k$-space, and there are at least two symps with rank at least $k+1$, but symps of rank at most $k$ are allowed. Fix $k\geq 2$. Then the Cartesian product of an arbitrary number of $k$-lacunary parapolar spaces, polar spaces and linear spaces, is a weakly $k$-lacunary parapolar space. If the direct factors are strong, then so is the product and, admittedly, it appears as the point-residual of a weakly $(k+1)$-lacunary parapolar space, which also could be strong, etcetera. 

For example, $\mathsf{E_{6,1}}(\mathbb{K})$ is $2$-lacunary, and so the Cartesian product with a thick line is weakly $2$-lacunary; hence $\mathsf{E_{8,7}}(\mathbb{K})$, or any homomorphic image, is weakly $3$-lacunary. The same argument with $\mathsf{D_{5,5}}(\mathbb{K})$ leads to $\mathsf{E_{7,6}}(\K)$. In fact, by taking the direct product with a projective space over $\mathbb{K}$, the geometries $\mathsf{E}_{n,6}(*)$ and $\mathsf{E}_{n,7}(*)$, $n\geq 8$, can be seen as weakly $3$-lacunary parapolar spaces. Moreover, we can consider a multiple direct product of $\mathsf{E_{6,1}}(\K)$ and/or $\mathsf{D_{5,5}}(\K)$ with several thick lines and projective spaces, and then the number of weakly $3$-lacunary parapolar spaces arising from non-spherical, non-Euclidean and non-hyperbolic Tits-buildings with a connected Coxeter diagram becomes  insurmountable. This motivates our present definition of $k$-lacunarity.

\section{Local recognition lemmas}\label{localrecognition}

In this section, we collect and prove some local recognition lemmas. There are two kinds of possible strategies that we can use. One of them is to use directly Tits' local approach theorem \cite{Tits:83}. We only prove the most involved case in detail. The other strategy is to use known characterization results and show that their assumptions are satisfied only knowing the local structure. The two methods together provide local recognition theorems for all parapolar spaces in the conclusion of our theorems, except for some with lacunary index $0$ (and there we will use Theorem~\ref{KasikovaShult2} and Proposition~\ref{prop:imbrex}).  

\begin{lemma}\label{KillEn}
 Let $\Omega=(X,\cL)$ be a (locally connected) parapolar space with the property that each point-residual is a homomorphic image of $\mathsf{D}_{n,n}(\K)$, for $n\geq 5$. Then $\Omega$ is a homomorphic image of $\mathsf{E}_{n+1,1}(*)$. The homomorphism is an isomorphism to $\mathsf{E}_{n+1,1}(\K)$, if $n\leq 6$. 
\end{lemma}
\begin{proof}
Let $L\in\cL$ be arbitrary. Then our assumption implies that $\Omega_L$ is isomorphic to $\mathsf{A}_{n-1,2}(\K)$. Since $\mathsf{A}_{n-1,2}(\K)$ has two types of maximal singular subspaces|of dimensions $2$ and $n-2$|we see that $\Omega$ also has two types of maximal singular subspaces|of dimensions $4$ and $n$, respectively. Also, since in $\Omega_L$ every line is contained in a unique maximal singular subspace of each type, in $\Omega$ each singular $3$-space is the intersection of a unique maximal singular $4$-space and a unique maximal singular $n$-space. Also, since in $\Omega_L$ every non-maximal singular plane is contained in a unique symp, every non-maximal singular $4$-space of $\Omega$ is also contained in a unique symp. Hence the points, lines, planes, $3$-spaces and $4$-spaces of a maximal singular subspace $U$ of $\Omega$ with $\dim U=n$ correspond to unique points, lines, planes, maximal singular $4$-spaces and symps of $\Omega$.

A \emph{chamber} of $\Omega$ is a set of $n+1$ objects consisting of a maximal singular subspace $U$ of dimension $n$ of $\Omega$ and a maximal flag in $U$ (that is, a nested sequence of $n$ non-trivial subspaces of $U$). We define the type of an element of a chamber $C$ as $2$, if it is the maximal singular subspace $U$ of dimension $n$, type $1$ if it is a point, and type $j$ if it is a subspace of dimension $j-2$ of $U$.   Then a chamber $C$ of $\Omega$ consists of $n+1$ elements of different type.  By the discussion in the previous paragraph, the elements of types $1$ up to $6$ of a chamber are a point $p$, a  maximal singular subspace $U$ of dimension $n$ containing $p$, a line $L$ with $p\in L\subseteq U$, a plane $\pi$ with $L\subseteq\pi\subseteq U$, a maximal singular subspace $W$ of dimension $4$ with $\pi\subseteq W$ and $\dim(U\cap W)=3$, and a symp $\xi$ with $W\subseteq\xi$ and $U\cap\xi$ a generator of $\xi$. We refer to this as the \emph{truncated representation} of the chamber.

Let $U$ be a singular $n$-space of $C$ and let $L\in\cL$ be contained in $U$. Let $\alpha$ be an $i$-space of $U$ containing $L$, $i\geq 4$. Then $\alpha$ corresponds to a unique subspace of $\Omega_L$ isomorphic to $\mathsf{A}_{i-1,2}(\K)$ and which intersects the maximal singular subspace $U_L$ of dimension $n-2$ of $\Omega_L$ corresponding to $U$ in a singular $(i-2)$-space. Conversely, a subspace of $\Omega_L$ isomorphic to $\mathsf{A}_{i-1,2}(\K)$ intersecting $U_L$ in a singular $(i-2)$-space corresponds to a unique $i$-space of $U$ containing $L$. This bijective correspondence can most easily be seen through the (partial) identification of the diagram of $U$ and the one of $\Omega_L$. The upshot is that a chamber of $\Omega$ can also be seen in a unique way as a set consisting of a point $p$, a line $L\ni p$, and a ``chamber'' of $\Omega_L$, that is, a set consisting of a singular $(n-2)$-space $U_L$ of $\Omega_L$, a point of $\Omega_L$ in $U_L$, a maximal singular plane of $\Omega$ containing $p$ and a nested sequence of subspaces $\Delta_j$ of $\Omega_L$ isomorphic to $\mathsf{A}_{j,2}(\K)$, $3\leq j\leq n-2$ intersecting $U_L$ in a singular $(j-1)$-space.  We refer to this as the \emph{line-residual representation} of the chamber.

Let $C$ and $C'$ be two chambers. Then we say that $C$ and $C'$ are \emph{$i$-adjacent}, $1\leq i\leq n+1$, $i\neq 2$, if $C$ and $C'$ contain the same maximal singular $n$-space, and if they also have all other elements in common, except possibly for their type $i$ elements. If $i=2$, then $C$ and $C'$ have distinct maximal singular $n$-spaces $U$ and $U'$ with $U\cap U'$ a singular plane $\pi$, but agree in all other elements, that is, they have the same point $p\in \pi$ and the same line $L\subseteq\pi$ and agree in all other elements of $\Omega_L$ in their line-residual representations, except possibly for $U$ and $U'$.  In symbols, we write $C\sim_iC'$.

We will now use the notion of a \emph{chamber system (of type $M$)}. Since this is only needed in the current proof, we refer the reader to \cite{Tits:83} for all definitions and background.  
Let $\cC$ be the set of chambers. Set $I=\{1,2,\ldots,n+1\}$. Then $\Gamma=(\cC,(\sim_i)_{i\in I})$ is a chamber system over $I$. We now show that $\Gamma$ is a chamber system of type $M$ (see \cite{Tits:83}), 
more exactly belonging to the diagram $\mathsf{E}_{n+1}$ as shown in the diagram of Figure~\ref{en+1}. According to Section~3.2 of \cite{Tits:83} it suffices to show that each residue of rank 2 is the chamber system of a (proper) projective plane or a generalized digon, and that $\Gamma$ is connected, that is, the graph $G_\Gamma$ with vertex set $\cC$ and adjacency relation the union over $I$ of the $i$-adjacency relations, is connected. 

For convenience we recall that the $J$-residue of a chamber $C\in\cC$, with $J\subseteq I$, is the connected component of $C$ in the graph with vertex set $\cC$ and edges the pairs of distinct $j$-adjacent chambers, with $j\in J$. A rank $k$ residue is a $J$-residue with $|J|=k$. Each $J$-residue is a (connected) chamber system over the type set $J$, in the obvious way. 

\begin{figure}[ht]
\begin{center}
\begin{tikzpicture}[scale=0.4]
\node at (0,0.3) {};
\node [inner sep=0.8pt,outer sep=0.8pt] at (-6,0) (1) {$\bullet$};
\node [inner sep=0.8pt,outer sep=0.8pt] at (-6,0.6) (1) {\footnotesize 1};

\node [inner sep=0.8pt,outer sep=0.8pt] at (-3,0) (3) {$\bullet$};
\node [inner sep=0.8pt,outer sep=0.8pt] at (-3,0.6) (3) {\footnotesize 3};

\node [inner sep=0.8pt,outer sep=0.8pt] at (0,0) (4) {$\bullet$};
\node [inner sep=0.8pt,outer sep=0.8pt] at (0,0.6) (4) {\footnotesize 4};

\node [inner sep=0.8pt,outer sep=0.8pt] at (3,0) (5) {$\bullet$};
\node [inner sep=0.8pt,outer sep=0.8pt] at (3,0.6) (5) {\footnotesize 5};

\node [inner sep=0.8pt,outer sep=0.8pt] at (6,0) (6) {$\bullet$};
\node [inner sep=0.8pt,outer sep=0.8pt] at (6,0.6) (6) {\footnotesize 6};

\node [inner sep=0.8pt,outer sep=0.8pt] at (9,0) (7) {$\bullet$};
\node [inner sep=0.8pt,outer sep=0.8pt] at (9,0.6) (7) {\footnotesize 7};

\node [inner sep=0.8pt,outer sep=0.8pt] at (0,-2) (2) {$\bullet$};
\node [inner sep=0.8pt,outer sep=0.8pt] at (0,-2.6) (2) {\footnotesize 2};

\node [inner sep=0.8pt,outer sep=0.8pt] at (12,0)  {$\bullet$};
\node [inner sep=0.8pt,outer sep=0.8pt] at (12,0.6) {\footnotesize $n$};

\node [inner sep=0.8pt,outer sep=0.8pt] at (15,0) {$\bullet$};
\node [inner sep=0.8pt,outer sep=0.8pt] at (15,0.6) {\footnotesize $n\!\!+\!\!1$};
\draw (-6,0)--(9,0);
\draw (0,0)--(0,-2);
\draw (12,0)--(15,0);
\draw (9,0)--(10,0);
\draw [dotted] (10,0)--(11,0);
\draw (11,0)--(12,0);
\end{tikzpicture}
\end{center}
\caption{The Coxeter diagram of type $\mathsf{E}_{n+1}$ with (extended) Bourbaki labeling\label{en+1}}
\end{figure}
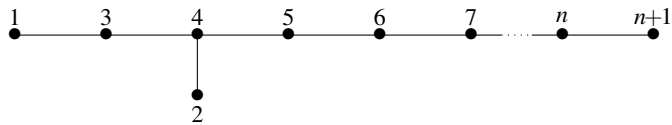

\textbf{Claim 1.} \emph{Let $i,j\in I$, $i\neq j$. Each $\{i,j\}$-residue of $\Gamma$ is the chamber system of a (proper) projective plane, if $i$ and $j$ are connected by an edge in Figure~\ref{en+1},  or a generalized digon otherwise.}

Suppose first that $2\notin\{i,j\}$. Then each residue of type $\{i,j\}$, say determined by the chamber $C$ containing the singular $n$-space $U$, is a residue of the chamber system of the projective space $U$ of dimension $n$ with diagram the one induced on $\{1,3,4,\ldots,n+1\}$ of $\mathsf{E}_{n+1}$ (see Figure~\ref{en+1}). This shows the claim if $2\notin\{i,j\}$. Now suppose $j=2$. If $i\notin\{1,3\}$, then similarly, the $\{2,i\}$-residue of the chamber $C$, say containing the line $L\in\cL$, is a residue of the chamber system of the projective space underlying the line Grassmannian $\Omega_L$, with diagram the one induced on $\{2,4,5,\ldots,n+1\}$ of $\mathsf{E}_{n+1}$ (see once again Figure~\ref{en+1}). If $i\in\{1,3\}$ and $C$ contains the maximal singular $4$-space $W$ as element of type $5$ then, using the truncated representation of $C$, its $\{2,i\}$-residue is a residue of the chamber system of the projective space $W$ with diagram the one induced on $\{1,3,4,2\}$ of $\mathsf{E}_{n+1}$. This shows the claim completely. 

\textbf{Claim 2.} \emph{The chamber system $\Gamma$ is connected.}

Indeed, let $C$ and $C'$ be two chambers. By connectivity of $\Omega$, we may assume that the respective elements $p$ and $p'$ of type $1$ of $C$ and $C'$ are collinear. Let $L''\in\cL$ containing both $p$ and $p'$. Let $C''$ be any chamber containing $L''$. It suffices to show that $C$ and $C''$ belong to the same connected component of $G_\Gamma$ (then similarly, also $C'$ and $C''$ belong to the same connected component). Clearly we may assume that $p\in C''$. Set $\{L\}= C\cap\cL$. Since $\Omega_p$ is connected, we may assume that $L$ and $L''$ belong to the same singular plane $\pi$. Since $\Omega_L$ is isomorphic to $\mathsf{A}_{n-1,2}(\K)$, the line-residual representations of the chambers imply that in the $(I\setminus\{1,3\})$-residue of $C$, the latter is connected to a chamber containing $\pi$. Similar for $C''$, and so we may assume that both $C$ and $C''$ contain $\pi$. Then $C''$ is $3$-adjacent to a chamber containing $\{p,L,\pi\}$, and hence we may assume that $L''=L$. Now the claim follows from the observation that the chamber system over $I\setminus\{1,3\}$ induced on the set of all chambers containing $\{p,L\}$ is isomorphic to the chamber system naturally associated to $\mathsf{A}_{n-1,2}(\K)$, and hence is connected.

\textbf{Claim 3.} \emph{Each rank $3$ residue of $\Gamma$ is the chamber system of a spherical building (whose type can be read off the $\mathsf{E}_{n+1}$ diagram of} Figure~\ref{en+1}). 

Similarly as in the proof of Claim 1, this follows for $\{i,j,k\}$-residues if $2\notin\{i,j,k\}$. Using the line-residual representation of chambers, Claim 3 follows for $\{2,i,j\}$-residues with $\{i,j\}\cap\{1,3\}=\emptyset$. The truncated representation of chambers implies Claim 3 for $\{2,i,j\}$-residues with $\{i,j\}\subseteq\{1,3,4,5\}$. Finally, suppose $i\in\{1,3\}$ and $j\in\{6,7,\ldots,n+1\}$. Then Claim 1 implies that each $\{1,i,j\}$-residue is the chamber system of a building of type $\mathsf{A_1\times A_1\times A_1}$. Claim 3 is proved.

Now Corollary~3 of \cite{Tits:83} implies that $\Gamma$ is the quotient by a group $G$ of the chamber system of a building of type $\mathsf{E}_{n+1}$. If $n=5$, then by \cite{Bro-Coh:83}, $G$ is trivial. If $n\geq 6$, then the diagram $\mathsf{E}_{n+1}$ does not admit non-trivial symmetries, hence the group $G$ is type-preserving, and so $G$ acts on the elements of respective types $1$ and $3$, implying that $\Omega$ is a homomorphic image of the $1$-Grassmannian of a building of type $\mathsf{E}_{n+1}$. 
Now let $n=6$. Then $G$ must also be trivial as the parapolar space $\mathsf{E_{7,1}}(\K)$ has diameter 3 and hence does not admit quotients that are parapolar spaces. 
\end{proof}

\begin{remark}
Lemma~\ref{KillEn} is proved in the course of the proof of Corollary~16.7.1 in \cite{Shu:12}, using the theory of truncated geometries and sheaves developed in Chapter~11 of \cite{Shu:12}. However, since Shult himself describes his proof as ``informal'' on page~582, and since the original paper \cite{Shu:05} mentioning Corollary~16.7.1 of \cite{Shu:12}  does not provide any proof, except for referring, without any explicit reference, to ``the theory of locally truncated geometries due to Ronan/Brouwer-Cohen'', we chose to provide a fully detailed proof, directly using Tits' results of \cite{Tits:83} rather than the detour via the theory of truncated geometries and sheaves, which is also based on \cite{Tits:83}.   
\end{remark}

Similarly, one shows the following recognition lemmas, the proofs of which are completely similar to, but simpler than, the proof of Lemma~\ref{KillEn}.

\begin{lemma}\label{KillDn}
Let $\Omega=(X,\cL)$ be a (locally connected) parapolar space with the property that each point-residual is isomorphic to $\mathsf{A}_{n,2}(\L)$, for $n\geq 4$ and some skew field $\L$. Then $\Omega$ is a homomorphic image of $\mathsf{D}_{n+1,n+1}(\K)$, with $\K=\L$ a field. The homomorphism is an isomorphism to $\mathsf{D}_{n+1,n+1}(\K)$, if $n\leq 8$. 
\end{lemma}

Lemma~\ref{KillDn} recognizes the point-residuals of the parapolar spaces we recognize in Lemma~\ref{KillEn}. Similarly, we can also recognize the point-residuals of $\mathsf{D}_{n,n}(\K)$ locally, and so starting with these point-residuals, we have three times in a row well determined extensions.
\begin{lemma}\label{KillAn2}
Let $\Omega=(X,\cL)$ be a (locally connected) parapolar space with the property that each point-residual is isomorphic to the Cartesian product of a thick line with $\mathsf{A}_{n,1}(\L)$, for $n\geq 2$ and some skew field $\L$. Then $\Omega$ is isomorphic to the building Grassmannian $\mathsf{A}_{n+2,2}(\L)$.  
\end{lemma}
    
\begin{lemma}\label{KillEnn}
Let $\Omega=(X,\cL)$ be a (locally connected) parapolar space with the property that each point-residual is isomorphic to $\mathsf{E}_{n,n}(\K)$, for $n\in\{6,7\}$ and some field $\K$. Then $\Omega$ is isomorphic to $\mathsf{E}_{n+1,n+1}(\K)$. 
\end{lemma}


We also want to recognize locally all long root geometries and close relatives that are parapolar spaces and have symplectic rank at least 3. To do this, it is convenient to use a result by Kasikova and Shult \cite{KasikovaShult}. Note that we already recognized $\mathsf{E_{7,1}}(\K)$ in Lemma~\ref{KillEn} and $\mathsf{E_{8,8}}(\K)$ in Lemma~\ref{KillEnn}. First we treat the classical cases.

\begin{lemma}\label{KillBn2}
Let $\Omega=(X,\cL)$ be a (locally connected) parapolar space with the property that each point-residual is isomorphic to the direct product of a thick line and a polar space (the isomorphism type of which may vary from point to point). Then $\Omega$ is isomorphic to $\mathsf{B}_{n,2}(*)$ or $\mathsf{D}_{n,2}(\K)$, $n\geq 4$, for some field $\K$.
\end{lemma}

\begin{proof}
We use \cite[Theorem 1]{KasikovaShult}. It suffices to prove the following three properties.

\textbf{Property 1.} \emph{ Given a point $q\in X$ not contained in a symp $\xi$, the singular subspace $q^\perp \cap \xi$ is never just a point.}\\
Indeed, suppose for a contradiction that $q^\perp\cap\xi$ is a point $p$. Then in $\Omega_p$, we find a point $p'$ and a symp $\xi'$ such that $p'^\perp\cap \xi'=\emptyset$. However, $\Omega_p$ is isomorphic to $L\times \Delta$, where $L$ is a thick line and $\Delta=(X',\cL')$ is a polar space. A point $p'=(x,x')$ of $\Omega_p$ is collinear to the point $(x'',x')$ of the symp $\{x''\}\times X'$, $x''\in L$, and to a point $(x,x''')$, with $x'''\in x'^\perp\cap L'$ in $\Delta$, of the symp $L\times L'$, with $L'\in\cL'$. Hence Property 1 is checked.

\textbf{Property 2.} \emph{Given a singular plane $\pi$ and line $M$ meeting $\pi$ at a point $p$, either $(i)$ every line of $\pi$ on $p$ lies in a common symp with $M$, or else $(ii)$ exactly one such line has this property.}\\
This translates to $\Omega_p$ as follows: For each point $p'$ and each line $L'$ not containing $p'$, either every point on $L'$ is at distance $2$ from $p'$, or exactly one is. Using the same notation as in the proof of Property 1, we set $p'=(x,x')$, with $x\in L$ and $x'\in X'$. Suppose first that $L'=\{y\}\times K'$, $K'\in\cL'$, with either $x\neq y$ or $x'\notin K'$. If $x=y$ or $x'\in K'$ then $p'$ and $L'$ are contained in a common symp. If $x\neq y$ and $x'\notin K'$, then $x'$ is at distance 2 from the points $(y,x'')$, with $x''\perp x'$ and $x''\in K'$. Secondly, suppose $L'=L\times\{x''\}$, $x''\in X'$, then $p'$ is only at distance 2 from $(x, x'')$ if $x'\notin x^\perp$, and at distance at most 2 from every point of $L'$ otherwise. This shows Property 2.

\textbf{Property 3.} \emph{Given any line $L$ on a point $p$, there exists at least one further line $N$ on $p$ such that $L^\perp \cap N^\perp = \{p\}$.}\\
This translates to $\Omega_p$ as follows: no point is at distance at most $2$ from every other point, which is clearly true. 

The assertion now follows from Theorem 1 of \cite{KasikovaShult}. 
\end{proof}

The following lemmas have similar proofs, and we leave them for the reader.

\begin{lemma}\label{KillF4}
Let $\Omega=(X,\cL)$ be a (locally connected) parapolar space with the property that each point-residual is isomorphic to a dual polar space of rank $3$. Then $\Omega\cong\mathsf{F_{4,1}}(*)$.
\end{lemma}


\begin{lemma}\label{KillE62}
Let $\Omega=(X,\cL)$ be a (locally connected) parapolar space with the property that each point-residual is isomorphic to the plane Grassmannian $\mathsf{A_{5,3}}(\L)$, for some skew field $\L$. Then $\L$ is a field and $\Omega\cong\mathsf{E_{6,2}}(\L)$.
\end{lemma}


\section{Parapolar spaces with lacunary index $0$}\label{sec0lac}
In this section $\Omega=(X,\cL)$ is a parapolar space of symplectic rank at least $d\geq 2$ and with lacunary index~$0$. Note that we do not assume local connectivity in this section, but we shall prove it. 

\subsection{Basic properties}
\label{basiczero}

\begin{lemma}\label{zerolocallyconnected}
Let $\Omega=(X,\cL)$ be a parapolar space of symplectic rank at least $d\geq 3$ with lacunary index $0$. Then for every point $p\in X$,  the point-residual $\Omega_p$ is a strong parapolar space of lacunary index $-1$. In particular, $\Omega$ is locally connected. 
\end{lemma}

\begin{proof}
By Axiom~(PPS1), there exists a point $p\in X$ contained in at least two symps. By Fact~\ref{pt-pps} and the absence of symps of rank 2, every connected component of $\Omega_p$ is a strong parapolar space, obviously with lacunary index $-1$. Since any two symps through $p$ have rank at least 3 and share a line, there is only one such component. Then by column $k=-1$ of Table~\ref{lacmin1}, proved in \cite{MinusOnePaper}, every line through $p$ is contained in at least two symps. Hence every point of $\Omega$ collinear to $p$ is contained in at least two symps and we can interchange its role with that of $p$. A connectivity argument now shows that $\Omega_x$ is a strong parapolar space with lacunary index $-1$, for every $x\in X$. Fact~\ref{conn} completes the proof of the lemma.
\end{proof}

\begin{lemma}\label{index0strong}
Parapolar spaces with lacunary index $0$ are strong.
\end{lemma}
\begin{proof}
Recall that in the case of symplectic rank 2 this is an assumption. When the symplectic rank is at least 3 this follows from Lemma~\ref{zerolocallyconnected} 
and the fact that strong parapolar spaces of lacunary index $-1$ all have diameter 2, see Table~\ref{lacmin1}. 
\end{proof}

We now bound the diameter.
The following lemma is Exercise 13.26 in  \cite{Shu:12}. We provide a proof for completeness. 

 \begin{lemma}\label{strong0diam3}
 Let $\Omega=(X,\cL)$ be a parapolar space with lacunary index $0$. Then $\diam  \Omega\leq 3$. 
 \end{lemma}

\begin{proof}
Suppose for a contradiction that $p_0\perp p_1\perp p_2\perp p_3\perp p_4$ are points with $\delta(p_0,p_4)=4$. The symps $\xi(p_0,p_2)$ and $\xi(p_2,p_4)$ (which really are symps by Lemma \ref{index0strong}) have $p_2$ and hence a line $L$ in common. In the symp $\xi(p_2,p_i)$ there is a point $q_i$ on $L$ collinear to $p_i$, $i=0,4$. Then we have $p_0\perp q_0\perp q_4\perp p_4$ and $\delta(p_0,p_4)\leq 3$, a contradiction. This proves the lemma.
\end{proof}

We now consider diameters 2 and 3 separately.

\subsection{Diameter 2 and minimum symplectic rank 2}\label{subsub22}

Let $\Omega=(X,\cL)$ be a (strong) parapolar space with lacunary index $0$, minimum symplectic rank 2, and diameter $2$. We first prove that the symplectic rank is uniformly 2.

\begin{lemma}
Let $\Omega=(X,\cL)$ be a (strong) parapolar space with lacunary index $0$, diameter $2$ and minimum symplectic rank $2$. Then $\Omega$ has uniform symplectic rank $2$.
\end{lemma}
\begin{proof}
Suppose for a contradiction that there is a symp $\xi$ of rank at least 3. By connectivity, we may assume that some point $p\in\xi$ is also contained in a symp of rank 2. We first claim that every line $L$ through $p$ is contained in a symp of rank 3. Indeed, let $x\in L\setminus\{p\}$ and pick a point $q$ in $\xi$ not collinear to $p$. Then the symp $\xi(x,q)$ shares a line $K$ with $\xi$ and hence $K$ contains a point $y$ collinear to $x$. Since $K$ does not contain $p$, we obtain a line $yp$ all of whose points are collinear to $x$. The claim now follows from Fact~\ref{planelinesymp}.

We next claim that $\Omega_p$ is a parapolar space. Indeed, our first claim implies that, if $L,M$ are lines through $p$, then they belong to respective symps of rank at least $3$, and these share a line by $0$-lacunarity. Consequently $\Omega_p$ is a connected point-line geometry. Hence Fact~\ref{pt-pps} implies that $\Omega_p$ is either a single line, or the set of lines of a single symp of rank at least $3$, or a parapolar space. Clearly the first two possibilities are impossible since there are at least two symps trough $p$ by assumption (one of rank 2 and one of rank at least 3). 
The claim is proved. 

Clearly, $\Omega_p$ is strong and has lacunary index $-1$. As can be seen in Table~\ref{lacmin1},  these all have diameter 2. Therefore every pair of lines of $\Omega$ through $p$ is contained in a symp of rank at least $3$, contradicting the existence of a symp of rank 2 through $p$. This contradiction shows that all symps through $p$ have rank 2. This proves the lemma.
\end{proof}

Hence, if $\Omega=(X,\cL)$ is a parapolar space with lacunary index $0$, minimum symplectic rank 2, and diameter $2$, then the symplectic rank is 2 and we clearly have an imbrex geometry. So we can apply Proposition~\ref{prop:imbrex}.

\begin{itemize}
\item
If the symps are thick, then we obtain $(\mathsf{GQ})$ of Main Result~\ref{main2}. An example that this really occurs is given by taking for $(X,\Sigma)$ (with $\Sigma$ the family of maximal singular subspaces) the dual of any Hermitian quadrangle in a $4$-dimensional projective space, and for $\Xi$ the family of subquadrangles induced by the hyperplanes. 
\item
If the symps are non-thick, then $\Omega$ is the direct product space of two linear spaces $Y$ and $Z$. Suppose both contain at least two lines, say $L_Y,K_Y$ are lines of $Y$ and $L_Z,K_Z$ are lines of $Z$. Obviously, we may choose $L_Y$ and $K_Y$ intersecting (say in the point $p_Y$) , and likewise $L_Z$ and $K_Z$ may assumed to have a point, say $p_Z$, in common with each other. Then the symps $L_Y\times L_Z$ and $K_Y\times K_Z$ intersect only at the point $(p_Y,p_Z)$, which is impossible by $0$-lacunarity. Hence one of $Y,Z$ is trivial and isomorphic to a thick line. 
\end{itemize}
This takes care of the white cells in the column $k=0$ of Table~\ref{lac0}.  
\subsection{Diameter 2 and symplectic rank at least 3}\label{subsub23}

In this subsection we show that parapolar spaces of diameter 2, symplectic rank at least 3 and lacunary index 0 are isomorphic to either $\mathsf{A_{2,4}}(\L)$, or $\mathsf{D_{5,5}}(\K)$, for some skew field $L$ and a field $\K$.  

\begin{lemma}\label{CCverified}
Let $\Omega=(X,\cL)$ be a parapolar space of diameter $2$, symplectic rank at least $3$ and lacunary index $0$. Then either all point-residuals are isomorphic to $\mathsf{A_{1.1}}(*)\times\mathsf{A_{2,1}}(*)$, or all point-residuals are isomorphic to $\mathsf{A_{4,2}}(\K)$. 
\end{lemma}
\begin{proof}
By Lemma~\ref{zerolocallyconnected}, the point-residual $\Omega_p$, for any $p\in X$, is a strong parapolar space with lacunary index $-1$. By the main result of \cite{MinusOnePaper}, see also Column $k=-1$ of Table~\ref{lacmin1}, it is one of $\mathsf{A_{1,1}}(\mathbb{*})\times \mathsf{A_{2,1}}(\mathbb{*}),\:\mathsf{A_{2,1}}(\mathbb{*})\times \mathsf{A_{2,1}}(\mathbb{*}),\: \mathsf{A_{4,2}}(\mathbb{L}),\: \mathsf{A_{5,2}}(\mathbb{L}), \: \mathsf{E_{6,1}}(\mathbb{K})$. Clearly, if $p\perp q$ in $\Omega$, then the parameters (singular rank, symplectic rank; see Table~\ref{lacmin1}) of $\Omega_p$ and $\Omega_q$ coincide, which implies, given the above list, that $\Omega_p$ and $\Omega_q$ are isomorphic. By connectivity we conclude that all point-residuals are isomorphic. 

Suppose $\Omega_p$ is one of $\mathsf{A_{2,1}}(\mathbb{*})\times \mathsf{A_{2,1}}(\mathbb{*}),\: \mathsf{A_{4,2}}(\mathbb{L}),\: \mathsf{A_{5,2}}(\mathbb{L}), \: \mathsf{E_{6,1}}(\mathbb{K})$ and note that these all have the property that some point is not collinear to any point of some symp. In other words, we can find a symp $\xi$ containing $p$ and a point $x$ with $x^\perp\cap\xi=\{p\}$.   Pick $q\in\xi\setminus p^\perp$ and let $\zeta$ be the symp containing $x$ and $q$. By $0$-lacunarity, $\zeta\cap\xi$ contains a line $L$, and hence $x$ is collinear to some point of $L$, which is different from $p$, a contradiction. 
The lemma is proved.
\end{proof}
The next proposition takes care of the white cells in the column $k=0$ of Table~\ref{lacmin1}.  
\begin{prop}
A parapolar space of diameter $2$, symplectic rank at least $3$ and lacunary index $0$ is isomorphic to either $\mathsf{A_{4,2}}(\L)$, or $\mathsf{D_{5,5}}(\K)$, for some skew field $L$ or a field $\K$.  
\end{prop}

\begin{proof}
This follows directly from Lemma~\ref{CCverified} and the recognition Lemmas~\ref{KillDn} and~\ref{KillAn2}.
\end{proof}

\subsection{Diameter 3}\label{subdiam3}

In this case, we verify that the conditions of Theorem~\ref{KasikovaShult2} are fulfilled. 

\begin{lemma}\label{nonempty}
Let $\Omega=(X,\cL)$ be a (strong) parapolar space with lacunary index $0$ and diameter $3$. Then for every point-symplecton pair $(x,\xi)$, we have $x^\perp \cap \xi \neq \emptyset$. In particular, the distance between a point and a line in $\Omega$ is at most $2$. 
\end{lemma}
\begin{proof}
Consider a point $q$ not in $\xi$, for which there is a path of length two, say $q\perp r\perp s$, for $r,s\in X$, with $s \in \xi$. We may assume $q\notin s^\perp$. By assumption the symps $\xi(q,s)$ and $\xi$ intersect in a line $L$. But then $q$ is collinear with a point on $L$, which also lies in $\xi$.

So we can keep shortening the path from $x$ to $\xi$, which exists by connectivity, and hence we have proved the first part of the lemma. For the second statement, let $x\in X$ and $L\in\cL$, then we include $L$ in symp $\xi$, obtain a point $q\in\xi$ with $q\perp x$, and so any point on $L$ collinear to $q$ is at distance at most $2$ from $x$.  
\end{proof}

\begin{lemma}\label{subspace}
Let $\Omega=(X,\cL)$ be a (strong) parapolar space with lacunary index $0$ and diameter $3$. Then the points at distance at most $2$ from a given point $p$ form a subspace, which is a geometric hyperplane.
\end{lemma}
\begin{proof}
Let $L$ be a line containing at least two distinct points at distance at most $2$ from $p$, say $x,y$. If $\delta(p,L)=1$ there is nothing to prove, so we may assume $\delta(p,x)=\delta(p,y)=2$. Since $\Omega$ is strong, there are symps $\xi(p,x)$ and $\xi(p,y)$, which by $0$-lacunarity intersect in a line $L_p$. 
If the points $x_p,y_p\in L_p$ collinear to respectively $x,y$, coincide, then $x\perp y\perp x_p\perp x$, so $x_p$ is collinear to $L$, i.e., $d(p,r)\leq 2$ for each point $r\in L$. If $x_p\neq y_p$, then $p$ and $L$ are contained in the symp $\xi(y,x_p)=\xi(x,y_p)$ and hence $\delta(p,L)=1$, concluding the proof that the points at distance at most 2 from $p$ form a subspace. By Lemma \ref{nonempty} and Proposition 11.5.21 of \cite{BuCo} this subspace is a geometric hyperplane.
\end{proof}

\begin{lemma}\label{finiterank}
Let $\Omega=(X,\cL)$ be a (strong) parapolar space with lacunary index $0$ and symplectic rank at least $3$. Then the maximal singular subspaces have finite dimension. 
\end{lemma}
\begin{proof}
The point-residuals of $\Omega$ are (strong) parapolar spaces of lacunary index $-1$, by Lemma~\ref{zerolocallyconnected}.  These all have maximal singular subspaces of finite projective dimension by the classification, see Column $k=-1$ of Table~\ref{lacmin1}. Whence the lemma.
\end{proof}  

\begin{prop}\label{prop0}
A parapolar space with lacunary index $0$ and diameter $3$ is one of the following:
\begin{compactenum}[$\bullet$]
\item $\mathsf{D_{6,6}}(\mathbb{K}), \mathsf{A_{5,3}}(\mathbb{L})$ or $\mathsf{E_{7,7}}(\mathbb{K})$;
\item a dual polar space of rank $3$ (that is, $\mathsf{B_{3,3}}(*)$);
\item a product geometry $L\times \Delta$, where $L$ is a thick line, and $\Delta$ is a  polar space of rank at least $2$. 
\end{compactenum}
\end{prop}

\begin{proof}
This follows from Theorem~\ref{KasikovaShult2} and Lemmas~\ref{nonempty},~\ref{subspace}, and~\ref{finiterank}.\end{proof}
This takes care of the grey cells in the columns $k=0$ of Tables~\ref{lacmin1} and~\ref{lac0}.  

\section{Parapolar spaces with lacunary index 1}\label{sec1lac}

Let $\Omega=(X,\cL)$ be a locally connected parapolar space with lacunary index $1$, strong if there are symps of rank 2. We start by showing that the symplectic rank is at least 3 (however, see Theorem~\ref{1lacnonstrong} for the non-strong case when the symplectic rank is 2). Then, since the point-residuals are $0$-lacunary parapolar spaces, we know them by Proposition~\ref{prop0}. In order to be able to apply the local recognition lemmas of Section~\ref{localrecognition}, we show that all such residuals are of the same type. 

\begin{lemma}\label{1,minrank>2}
Let $\Omega=(X,\cL)$ be a strong parapolar space of minimum symplectic rank $2$. Then $\Omega$ is not $1$-lacunary.
\end{lemma}

\begin{proof}
Suppose for a contradiction that $\Omega$ is $1$-lacunary. Let $\xi$ be any symp of $\Omega$ of rank 2. By connectivity, there is a symp $\xi'$ intersecting $\xi$ non-trivially. Then $\xi \cap \xi'$ is a point $p$, as otherwise $1$-lacunarity forces $\xi \cap \xi'$ to contain a plane, which is impossible as $\xi$ has rank 2. Let $L$ and $L'$ be lines through $p$ in $\xi$ and $\xi'$, respectively. Since $\Omega$ is strong, there is a symp through $L$ and $L'$, which intersects $\xi$ in $L$, contradicting what we have just deduced.
\end{proof}
Taken together with the previous section, and in view of the main results of \cite{MinusOnePaper}, this completes the proof of Main Result~\ref{main2}. It also shows Corollary~\ref{corMR}$(i)$, by an obvious inductive argument.



\begin{lemma}\label{1,uniformres}
Let $\Omega=(X,\cL)$ be a locally connected parapolar space with lacunary index $1$ and symplectic rank at least $3$. Then, for each two points $p,q\in X$, the point-residuals $\Omega_p$ and $\Omega_q$ are Lie incidence geometries of the same Coxeter type.  
\end{lemma}

\begin{proof}
Take any point $p \in X$ and consider the point-residual $\Omega_p$. Since $\Omega$ is locally connected, Fact~\ref{conn} implies that $\Omega_p$ is connected, and then we obtain from Fact~\ref{pt-pps} that $\Omega_p$ is a (strong) parapolar space. Clearly, $\Omega_p$ is 0-lacunary. Moreover, the singular subspaces of $\Omega_p$ are projective, since this is the case for $\Omega$ (cf.\ Fact~\ref{planelinesymp}). Hence $\Omega_p$ is as in Tables~\ref{lacmin1} and~\ref{lac0} (columns corresponding to $k=0$), except that the $\mathsf{GQ}$-case does not occur, and that the linear space ``$\mathsf{LS}$'' is a projective space. Consequently, $\Omega_p$ is a Lie incidence geometry. One also observes that we can distinguish the entries of these columns by their symplectic and singular ranks. 

Now let $q\in X$ be collinear with $p$. We claim that $\Omega_q$ has the same symplectic and singular ranks as $\Omega_p$. Indeed, each symp through the line $pq$ corresponds with a unique symp of $\Omega_p$ through $q$ and with a unique symp $\Omega_q$ through $p$, and vice versa. Hence there is a bijective correspondence between the symps of $\Omega_p$ through $q$ and of $\Omega_q$ through $p$; likewise for the maximal singular subspaces. Each point in $\Omega_p$ (resp.\ $\Omega_q$) plays the same role, so the local parameters determine the global ones, proving the claim. So $\Omega_p$ and $\Omega_q$ have the same Coxeter type indeed, and by connectivity, the lemma follows.\end{proof}

\begin{prop}\label{prop11}\label{prop12}
Let $\Omega=(X,\cL)$ be a locally connected parapolar space with lacunary index $1$ and symplectic rank at least $3$.  Then either 
\begin{compactenum}[$(i)$] \item $\Omega$ is one of the non-strong parapolar spaces $\mathsf{B}_{n,2}(*)$, $\mathsf{D}_{n,2}(\mathbb{K})$, $\mathsf{E_{6,2}}(\K)$, $\mathsf{E_{7,1}}(\K)$, $\mathsf{E_{8,8}}(\K)$ or $\mathsf{F_{4,1}}(*)$, $n\geq 4$, $\K$ any field, or
%
\item $\Omega$ is one of the strong parapolar spaces $\mathsf{A}_{n,2}(\mathbb{L})$, $\mathsf{D_{5,5}}(\mathbb{K})$ or $\mathsf{E_{6,1}}(\mathbb{K})$, $n\geq 4$, $\K$ any field.
\end{compactenum}
\end{prop}

\begin{proof}
By Lemma~\ref{1,uniformres}, all point-residuals have the same type, and are moreover $0$-lacunary. We can now use the classification of $0$-lacunary parapolar spaces again. Suppose first that the point-residuals have diameter $3$. 
From Tables~\ref{lacmin1} and~\ref{lac0}, we deduce the following possibilities ($p\in X$ arbitrary).
\begin{compactenum}[$\bullet$]
\item $\Omega_p$ is the Cartesian product of a thick line and a polar space. Then, by Lemma~\ref{KillBn2}, $\Omega$ is one of $\mathsf{B}_{n,2}(*)$ or $\mathsf{D}_{n,2}(\K)$, $n\geq 4$, $\K$ any field.
\item $\Omega_p$ is dual polar space of rank 3. Then, by Lemma~\ref{KillF4}, $\Omega\cong\mathsf{F_{4,1}}(*)$.
\item $\Omega_p$ is one of $\mathsf{A_{5,3}}(\L)$, $\mathsf{D_{6,6}}(\K)$ or $\mathsf{E_{7,7}}(\K)$, $\L$ a skew field, $\K$ a field. Then, by Lemmas~\ref{KillE62},~\ref{KillEn} and~\ref{KillEnn}, respectively, $\Omega$ is either $\mathsf{E_{6,2}}(\L)$, $\mathsf{E_{7,1}}(\K)$ or $\mathsf{E_{8,8}}(\K)$, where $\L$ is necessarily a field. 
\end{compactenum}
This shows $(i)$. Now suppose the point-residuals have diameter 2. Then Tables~\ref{lacmin1} and~\ref{lac0} give the following possibilities, noting that singular subspaces are projective (and again $p\in X$ is arbitrary):
\begin{compactenum}[$\bullet$]
\item $\Omega_p$ is the Cartesian product of a thick line and a projective space of dimension at least $2$. Then, by Lemma~\ref{KillAn2}, $\Omega\cong\mathsf{A}_{n,2}(\L)$, $n\geq 4$ and $\L$ some skew field.
\item $\Omega_p$ is one of $\mathsf{A_{4,2}}(\L)$ or $\mathsf{D_{5,5}}(\K)$, $\L$ a skew field, $\K$ a field. Then, by Lemmas~\ref{KillDn} and~\ref{KillEn}, respectively, $\Omega$ is either $\mathsf{D_{5,5}}(\K)$ or $\mathsf{E_{6,1}}(\K)$, respectively (with $\L$ necessarily a field). 
\end{compactenum} This yields $(ii)$ and completes the proof of the proposition.
\end{proof}
Proposition~\ref{prop11}$(i)$ takes care of the grey cells in the columns $k=1$ of Tables~\ref{lacmin1} and~\ref{lac0}. Also, Proposition~\ref{prop11}$(ii)$ takes care of the white cells in the columns $k=1$ of Tables~\ref{lacmin1} and~\ref{lac0}.




\section{Parapolar spaces with lacunary index at least 2}\label{sec>2}

Let $\Omega=(X,\cL)$ be a locally connected parapolar space with lacunary index $k$, $k\geq 2$. By definition, the symplectic rank is at least $k+1$. 
%
%
We can now conclude the proof of our Main Theorem, using an induction on the lacunary index $k$, starting with the already established case $k=1$.
\begin{theorem}
Let $\Omega=(X,\cL)$ be a locally connected parapolar space with lacunary index $k$, $k \geq 2$, and minimal symplectic rank $d\geq k+1$. Then $\Omega$ is one of the parapolar spaces mentioned in Tables~\emph{\ref{lacmin1}} and~\emph{\ref{lac0}} in the columns corresponding to $k=2,3,4,5$.
\end{theorem}

\begin{proof}
We prove it step by step starting with the case $k=2$. Consider any point-residual $\Omega_p$, $p\in X$. Then, by local connectivity, $\Omega_p$ is a strong parapolar space with lacunary index $1$ and minimal symplectic rank $d'\geq d-1\geq k=2$. By Lemma~\ref{1,minrank>2}, $d'\geq 3$, and so, since $\Omega_p$ cannot be one of the parapolar spaces mentioned in the conclusion of Proposition~\ref{prop11}, since these are all non-strong, $\Omega_p$ is one of $\mathsf{A}_{n,2}(\mathbb{L})$, $\mathsf{D_{5,5}}(\mathbb{K})$, or $\mathsf{E_{6,1}}(\mathbb{K})$. 

Suppose first that $\Omega_p\cong \mathsf{E_{6,1}}(\mathbb{K})$. Then every point $q$ collinear to $p$ is contained in a symplecton of rank $5$, and hence, since $\Omega_q$ also has to be one of the above parapolar spaces and $\mathsf{E_{6,1}}(\mathbb{K})$ is the only one with symplectic rank 5, we see that $\Omega_q\cong\mathsf{E_{6,1}}(\mathbb{K})$. By connectivity all point-residuals are isomorphic to $\mathsf{E_{6,1}}(\mathbb{K})$. Now Lemma~\ref{KillEnn} implies $\Omega\cong\mathsf{E_{7,7}}(\mathbb{K})$. 
Similarly, if $\Omega_p\cong\mathsf{D_{5,5}}(\K)$, then, using Lemma~\ref{KillEn}, we conclude $\Omega\cong\mathsf{E_{6,1}}(\K)$.  
This already yields the column $k=2$ of Table~\ref{lacmin1}. 

Now suppose $\Omega_p\cong\mathsf{A}_{n,2}(\mathbb{L})$, $n\geq 4$. Then the maximum singular subspaces containing $p$ have dimension $n-1$. Since every line through $p$ is contained in such a singular subspace, the same is true for every point collinear to $p$. By connectivity, we conclude that all point-residuals are isomorphic to $\mathsf{A}_{n,2}(\mathbb{L})$. Now Lemma~\ref{KillDn} implies that $\L$ is a field and $\Omega\cong\mathsf{D}_{n+1,n+1}(\mathbb{K})^h$. This yields the column $k=2$ of Table~\ref{lac0}.

Now suppose $k=3$. Pick $p\in X$. Then $\Omega_p$ is a connected $2$-lacunary parapolar space, and thus one of $\mathsf{D}_{n,n}(\mathbb{K})^h$, $n\geq 5$,  $\mathsf{E_{6,1}}(\mathbb{K})$, or $\mathsf{E_{7,7}}(\mathbb{K})$. Similar arguments as for the case $k=2$ lead in the two latter cases to $\Omega\in\{\mathsf{E_{7,7}}(\mathbb{K}),\mathsf{E_{8,8}}(\mathbb{K})\}$ (using Lemma~\ref{KillEnn}), that is, the column $k=3$ in Table~\ref{lacmin1}. In the former case, a similar argument, together with Lemma~\ref{KillEn}, yields $\Omega\cong\mathsf{E}_{n+1,1}(*)^h$, where $\Omega\cong\mathsf{E}_{n+1,1}(\K)$ for $n\in\{5,6\}$. This is column $k=3$ in Table~\ref{lac0}. 

Next, suppose $k=4$. Again pick $p\in X$. Then $\Omega_p$ is a connected $3$-lacunary strong parapolar space. It is easy to see that $\mathsf{E}_{n+1,1}(*)^h$ is not strong for $n\geq 6$, as the point-residual of a strong parapolar space has diameter $2$. Also, $\mathsf{E_{8,8}}(\K)$ is not strong as a long root geometry. Hence, this time, $\Omega_p$ is isomorphic to $\mathsf{E_{6,1}}(\K)$ or $\mathsf{E_{7,7}}(\K)$. Similarly as above, using Lemma~\ref{KillEnn}, we conclude $\Omega\in\{\mathsf{E_{7,7}}(\mathbb{K}),\mathsf{E_{8,8}}(\mathbb{K})\}$, that is, the column $k=4$ of Tables~\ref{lacmin1} and~\ref{lac0}. 

Similarly, the case $k=5$ leads to $\Omega\cong\mathsf{E_{8,8}}(\K)$, that is, the column $k=5$ of Table~\ref{lac0}. 

Since the only $5$-lacunary locally connected parapolar spaces are non-strong, there do not exist $k$-lacunary locally connected parapolar spaces with $k\geq 6$. 

This completes the proof of the theorem, and of Main Result~\ref{main1}.
\end{proof}

 \appendix
 
 \section{Appendix: Locally disconnected parapolar spaces}\label{ap}
 
 It is well known, and attributed to folklore by Shult in \cite{Shu:12}, that a gamma-space can be assembled in a unique way from locally connected ``components'', see \cite[Theorem~3.6.1]{Shu:12}. Shult then uses this result to comment on parapolar spaces which are not locally connected, see \cite[Section~13.5]{Shu:12}. The answer given there, in particular the discussion in Subsection~3.6.2, is perhaps slightly too vague and too general. 
 
 In this note we revise Shult's results directly in the theory of parapolar spaces, providing full detailed proofs, although, as Shult notes, a lot of people have independently rediscovered some parts of it, but, as Shult also mentions in a footnote \cite[p. 69]{Shu:12}, also the ``easy results'' deserve detailed proofs and clearly stated corollaries. Our aim is to revise Shult's principles in \cite[3.6.2]{Shu:12}, showing that they do not suffice in certain situations, and illustrate the need for this revision with a few examples. 

 \subsection{Unbuttoning of parapolar spaces}

Let $\Omega=(X,\cL)$ be an arbitrary parapolar space with symplectic rank at least 3.  For each point $p\in X$, we denote by $\mathfrak{C}_p$ the set of connected components of $\Omega_p$ (see Definition~\ref{point-residue} and Fact~\ref{pt-pps}). By definition, $\Omega$ is \emph{locally connected} if $\Omega_p$ is connected for all $p\in X$; otherwise $\Omega$ is \emph{locally disconnected}. The following construction introduces a copy of a point $p$ for each connected component of $\Omega_p$.

\begin{con}
Let $\Omega=(X,\cL)$ be an arbitrary parapolar space with symplectic rank at least 3. The \emph{unbuttoning of $\Omega$} is defined as the following point-line geometry $\widetilde{\Omega}=(\widetilde{X},\widetilde{\cL})$:
\begin{compactenum}[$\bullet$]
\item $\widetilde{X}=\{(p,\Upsilon):p\in X\mbox{ and }\Upsilon\in \mathfrak{C}_p\}$;
\item for each line $L\in\cL$, we define $\widetilde{L}=\{(p,\Upsilon)\in\widetilde{X}:p\in L\in\Upsilon\}$,
\item $\widetilde{\cL}=\{\widetilde{L}:L\in\cL\}$.
\end{compactenum}
\end{con}
So two points $(p_1,\Upsilon_1)$ and $(p_2,\Upsilon_2)$, with $\Upsilon_i \in \Omega_{p_i}$ for $i=1,2$, are collinear in $\widetilde{\Omega}$ if and only if $p_1 \perp p_2$ and the line $p_1p_2$ is an element of both $\Upsilon_1$ and $\Upsilon_2$. {For a point $(p,\Upsilon)$, we call $p$ the \emph{first coordinate of $(p,\Upsilon)$}.}

This procedure yields a disjoint union of locally connected (para)polar spaces:

\begin{prop}\label{unbutton}
Let $\Omega=(X,\cL)$ be a not necessarily locally connected parapolar space of symplectic rank at least $3$. Then its unbuttoning $\widetilde{\Omega}$ is the disjoint union of locally connected (para)polar spaces.  
\end{prop} 

\begin{proof}
We verify the axioms of a parapolar space except that in Axiom~(PPS1) we do not require that there is a point-line pair $(p,L)$ such that no point of $L$ is collinear to $p$, nor do we require that $\widetilde{\Omega}$ is connected, instead we will in the end consider its connected components.

\begin{itemize}
\item[(PPS1)] Suppose $(p,\Upsilon)\in\widetilde{X}$ and $\widetilde{L}\in\widetilde{\cL}$ are such that $(p,\Upsilon)\notin \widetilde{L}$ is collinear to at least two points of $\widetilde{L}$. Let $(x^*, \Upsilon^*)$ be any point of $\widetilde{L}$. In $\Omega$, at least two points of $L$ are collinear to $p$, so $\<p,L\>$ is a plane $\pi$. This means that each line of $\pi$ through $p$ belongs to $\Upsilon$ (in particular,  $px^* \in \Upsilon$) and likewise each line of $\pi$ through $x^*$ is contained in $\Upsilon^*$ (in particular, $px^* \in \Upsilon^*$). Consequently, $(p,\Upsilon)$ and $(x^*,\Upsilon^*)$ are contained in $\widetilde{px^*}$ and as such they are indeed collinear in $\widetilde{\Omega}$.


\item[(PPS2)]  
Let $(p_i,\Upsilon_i)\in\widetilde{X}$, $i=1,2$, be two non-collinear points of $\widetilde{\Omega}$ collinear to at least one common point $(x_1,\Sigma_1)$ of $\widetilde{\Omega}$. We claim that $p_1$ and $p_2$ are not collinear in $\Omega$. Indeed, suppose they are. Since $(p_1,\Upsilon_1)$ is collinear to $(x_1,\Sigma_1)$, the line $p_1x_1$ belongs to $\Upsilon_1$. As $x_1$ is collinear to $p_2$, the line $p_1p_2$ lies in $\Upsilon_1$ too. Likewise we obtain $p_1p_2 \in \Upsilon_2$. But then the points $(p_i,\Upsilon_i)$, $i=1,2$, both belong to $\widetilde{p_1p_2}$, a contradiction.  Our claim follows.
Now suppose that both $(p_i,\Upsilon_i)$, $i=1,2$, are collinear to a second point $(x_2,\Sigma_2)$, with $(x_1,\Sigma_1)\neq (x_2,\Sigma_2)$. Since $x_ip_1\in\Sigma_i$ for $i=1,2$ and $\Sigma_1 \cap \Sigma_2 = \{x_1\}$ if $x_1=x_2$, we deduce $x_1\neq x_2$.

We now show that the convex closure $C$ of $(p_1,\Upsilon_1)$ and $(p_2,\Upsilon_2)$ is a polar space canonically isomorphic to the symp $\xi:=\xi(p_1,p_2)$.  To that aim, we have to show two claims. 

 \textbf{Claim 1:} \emph{ If $x\in\xi(p_1,p_2)$ and if we denote by $\Sigma_{x,\xi}$ the component of $\Omega_x$ containing the lines of $\xi$ through $x$, then $(x,\Sigma_{x,\xi})$   belongs to $C$. }\\
Indeed, by Fact~\ref{convexclosure}, it suffices to show that that $(x,\Sigma_{x,\xi}) \in C$ for all points $x$ which are contained in a line joining $p_1$ or $p_2$ with a point of $p_1^\perp \cap p_2^\perp$. Suppose first that $x \in p_1^\perp \cap p_2^\perp$. Then, firstly, $xp_1$ and $xp_2$ belong to $\Sigma_{x,\xi}$ by our assumption on $\Sigma_{x,\xi}$. Secondly, $p_ix$ belongs to $\Upsilon_i$, $i=1,2$, because $p_ix$ lies in the same connected component of $\Omega_{p_i}$ as $p_ix_1$ and $p_ix_2$ (these lines all lie in $\xi$), $i=1,2$. This shows that the point $(x,\Sigma_{x,\xi})$ is collinear to $(p_i,\Upsilon_i)$, $i=1,2$, and hence belongs to $C$ indeed. Similarly one can now show that each point $x'$ on the line $p_ix$ is such that $(x',\Sigma_{x',\xi})$ is on the line joining $(p_i,\Upsilon_i)$ and $(x,\Sigma_{x,\xi})$, $i=1,2$, and hence $(x',\Sigma_{x',\xi})\in C$ too. This shows Claim 1.
\medskip

 \textbf{Claim 2:} \emph{If  $(y,\Upsilon)\in C$, then $y\in\xi$, and $\Upsilon$ is the component of $\Omega_y$ containing the lines of $\xi$ through~$y$.}\\
Indeed, let $C'$ denote the set of points $(x,\Sigma_{x,\xi})$, with $x\in\xi$. Let $\rho$ be the projection map $C'\rightarrow\xi:(x,\Sigma)\mapsto x$.  Then $\rho$ is an isomorphism of point-line geometries: the first paragraph implies that $\rho$ preserves collinearity and is injective; surjectivity follows by definition of $C'$. Hence $C'$ is a polar space containing $(p_1,\Upsilon_1)$ and $(p_2,\Upsilon_2)$ and therefore $C'=C$. 

This concludes the verification of Axiom~(PPS2).

\item[(PPS3)] Let $\widetilde{L}$ be a line of $\widetilde{\Omega}$.  Then $L\in\cL$ is contained in some symp $\xi$. We consider two points $p_1,p_2\in\xi$ at distance $2$ and with $p_1\in L$. We showed above that $(p_1,\Sigma_{p_1,\xi})$ and $(p_2,\Sigma_{p_2,\Sigma})$ determine a symp $\widetilde{\xi}$ in $\widetilde{\Omega}$, which contains precisely the points $(x,\Sigma_{x,\xi})$ with $x \in \xi$, so in particular those with $x \in L$. Since $L \subseteq \xi$, we have that $L \subseteq \Sigma_{x,\xi}$ for all $x\in L$, showing that $\widetilde{L}$ belongs to $\widetilde{\xi}$.
\end{itemize}

This shows that each connected component $\omega$ of $\widetilde{\Omega}$ is a (para)polar space. The fact that $\omega$ is locally connected follows immediately from the definition of $\widetilde{\Omega}$.   This proves the proposition.
\end{proof}

\begin{remark}
Note that, in Proposition~\ref{unbutton},  we restricted ourselves to parapolar spaces with symplectic rank at least 3, in contrast to the general theory of gamma spaces (see Section 3.6 of \cite{Shu:12}). The reason is that, although for gamma spaces this might make sense, for the more restricted class of parapolar spaces we are not aware of any sensible way of unbuttoning  so that both symps and singular subspaces can be reassembled afterwards.  
\end{remark} 

We are now interested in a reverse procedure. Which parapolar spaces can we obtain by collecting locally connected (para)polar spaces and identifying certain points? 

The following lemma is necessary to make the construction universal. Basically it says that, in $\Omega$, you cannot walk from a point $p$ to itself in less than five steps using two different components of $\Omega_p$ to start and come back in.

\begin{lemma}\label{universal}
Let $\Omega=(X,\cL)$ be a not necessarily locally connected parapolar space with symplectic rank at least $3$. Let $\widetilde{\Omega}$ be its unbuttoning. Let $p\in X$ be such that $\Omega_p$ is disconnected and let $\Upsilon^{(p)}_1$ and $\Upsilon^{(p)}_2$ be two distinct connected components of $\Omega_p$. Let $q,r,s\in X\setminus\{p\}$ be arbitrary (not necessarily distinct) and let $\Upsilon^{(q)}_1,\Upsilon^{(q)}_2,\Upsilon^{(r)}_1,\Upsilon^{(r)}_2,\Upsilon^{(s)}_1,\Upsilon^{(s)}_2$ be not necessarily distinct respective connected components (with self-explaining notation) of $\Omega_q,\Omega_r,\Omega_s$.  Then 
\[\ell:=\delta((p,\Upsilon^{(p)}_2),(q,\Upsilon^{(q)}_1))+\delta((q,\Upsilon^{(q)}_2),(r,\Upsilon^{(r)}_1))+\delta((r,\Upsilon^{(r)}_2),(s,\Upsilon^{(s)}_1))+\delta((s,\Upsilon^{(s)}_2),(p,\Upsilon^{(p)}_1))\geq5\] (points in different components of $\widetilde{\Omega}$  have distance $\delta=\infty$, which is by definition larger than any positive number).   
\end{lemma}

\begin{proof}
Suppose for a contradiction that $\ell\leq 4$. We examine the case $\ell=4$, leaving the easier cases $\ell=1,2,3$ to the interested reader. The assumption $\ell=4$ allows us to also assume that $$\delta((p,\Upsilon^{(p)}_2),(q,\Upsilon^{(q)}_1))=\delta((q,\Upsilon^{(q)}_2),(r,\Upsilon^{(r)}_1))=\delta((r,\Upsilon^{(r)}_2),(s,\Upsilon^{(s)}_1))=\delta((s,\Upsilon^{(s)}_2),(p,\Upsilon^{(p)}_1))=1,$$ since, if some of these distances would be 0, then another distance must be at least 2 and we can insert a chain of points consecutively at distance 1, rename, and get the above assumption back. 

By the definition of lines in $\widetilde{\Omega}$ we then obtain $p\perp\ q\perp r\perp s\perp p$. First note that the lines $pq$ and $ps$ belong to $\Upsilon_2^{(p)}$ and $\Upsilon_1^{(p)}$, respectively. By assumption, $\Upsilon_1^{(p)} \neq \Upsilon_2^{(p)}$. This already implies $q \neq s$. It also implies that $q$ cannot be collinear to $s$, for then $\<p,q,s\>$ would be a projective plane (cf.\ Fact~\ref{subspacesymp}), yielding $\Upsilon_1^{(p)} = \Upsilon_2^{(p)}$ after all. However, if $q$ and $s$ are not collinear, they determine a symp $\xi$ since $p \neq r$, clearly containing the lines $pq$ and $ps$, which again leads to $\Upsilon_1^{(p)} = \Upsilon_2^{(p)}$. This contradiction proves the lemma.
\end{proof}

\subsection{Buttoning of parapolar spaces}

\begin{con}\label{buttoned}
Let $\mathcal{F}=\{\Omega_i=(X_i,\cL_i):i\in I\}$ be a family of (disjoint) locally connected (para)polar spaces (of symplectic rank at least $2$), over some non-empty index set $I$.  
Let $\mathcal{R}$ be an equivalence relation on the union $\widetilde{X}=\bigcup_{i\in I}X_i$ of the sets of points of all members of $\mathcal{F}$, satisfying the following two conditions (C1) and (C2).
\begin{itemize}
\item[(C1)] Let $\widetilde{p},\widetilde{q},\widetilde{r},\widetilde{s}$ be four (not necessarily distinct, but $\widetilde{p} \notin\{\widetilde{q},\widetilde{r},\widetilde{s}\}$) equivalence classes with respect to $\mathcal{R}$, and let $p_1,p_2\in \widetilde{p}$, with $p_1\neq p_2$. If $q_1,q_2\in \widetilde{q}$, $r_1,r_2\in \widetilde{r}$ and $s_1,s_2\in \widetilde{s}$, then $$\delta(p_2,q_1)+\delta(q_2,r_1)+\delta(r_2,s_1)+\delta(s_2,p_1)\geq 5.$$
\item[(C2)] The graph with vertex set $\mathcal{F}$, where two vertices $\Omega_i$ and $\Omega_j$, $i,j\in I$, are adjacent if some point of $\Omega_i$ is contained in the same equivalence class as some point of $\Omega_j$, is connected.  
\end{itemize}
Set $X=\widetilde{X}/\mathcal{R}$. For each line $L$ contained in some member of $\mathcal{F}$, we put $\widetilde{L}:=\{ \widetilde{p} \mid p \in L\}$ and define $\cL$ as $\{\widetilde{L} \mid L \in \cL_i \text{ for some } i\in I\}$. 
Then we denote the geometry $\Omega=(X,\cL)$ by $\Omega(\mathcal{F},\mathcal{R})$. If $\mathcal{R}$ is non-trivial, then we call $\Omega$ a \emph{buttoned geometry}. The members of the family $\mathcal{F}$ are called \emph{sheets}. If $\mathcal{C}$ is a class of polar and parapolar spaces such that each member of $\mathcal{F}$ is contained in $\mathcal{C}$, then we say that $\Omega$ \emph{is buttoned over $\mathcal{C}$}, and $\mathcal{F}$ is called the \emph{sheet space} of $\Omega$ (note that members of $\mathcal{C}$ are allowed to occur multiple times in $\mathcal{F}$). {If $\mathcal{F}=\{\Omega^*\}$ for some parapolar space $\Omega^*$, then $\Omega$ is called \emph{self-buttoned}.} 
%
\hfill $\blacksquare$
\end{con}

\begin{lemma} \label{uniqueline} With the notation of \emph{Construction~\ref{buttoned}}, distinct lines $L$ and $L'$ define distinct sets $\widetilde{L}$ and $\widetilde{L'}$. \end{lemma}

\begin{proof} Indeed, suppose $\widetilde{L}=\widetilde{L'}$ and take points $\widetilde{p},\widetilde{q} \in \widetilde{L}=\widetilde{L'}$ with $\widetilde{p} \neq \widetilde{q}$. Let $p_1,q_2$ be (distinct) points on $L$ and $p_2,q_1$ (distinct) points on $L'$ with $p_1,p_2 \in \widetilde{p}$ and $q_1,q_2 \in \widetilde{q}$. Then Condition $(C1)$, with $\widetilde{r}=\widetilde{s}=\widetilde{q}$, $r_1=s_2=q_2$ and $r_2=s_1=q_1$  implies $p_1=p_2$. Since $\widetilde{p} \in \widetilde{L}$ was arbitrary, we obtain $L=L'$. The claim follows. \end{proof}

We have the following result. 
\begin{prop}\label{crucial}
{Let $\Omega=\Omega(\mathcal{F},\mathcal{R})$ be a buttoned geometry with sheet space $\mathcal{F}$ and non-trivial equivalence relation $\mathcal{R}$.}
Then $\Omega$ is a locally disconnected parapolar space. {Moreover, there is a natural bijection between the set of symplecta of $\Omega$ and the set of symplecta of all members of $\mathcal{F}$ and corresponding symps are isomorphic and hence have the same rank; in particular, if all members of $\mathcal{F}$ have (symplectic) rank at least 3, then so does $\Omega$.} 
\end{prop}

\begin{proof}
We verify the axioms of a parapolar space for $\Omega:=\Omega(\mathcal{F},\mathcal{R})$.
\begin{itemize}
\item[(PPS1)] Condition (C2) implies immediately that $\Omega$ is connected.  Let $\widetilde{L}\in\cL$ and $\widetilde{p}\in X$ with $\widetilde{p} \notin \widetilde{L}$ be such that $\widetilde{p}$ is collinear to at least two points $\widetilde{q},\widetilde{r}\in \widetilde{L}$. The definition of $\mathcal{L}$ yields (unique) points $p_1,p_2 \in \widetilde{p}$, $q_1,q_2 \in \widetilde{q}$, $r_1,r_2 \in \widetilde{r}$ with $p_1 \perp q_2$, $q_1 \perp r_2$, $r_1 \perp p_2$. By Condition~(C1), $p_1=p_2$, $q_1=q_2$ and $r_1=r_2$. It follows that $q_1r_1=L$ and hence $p_1$ is collinear to each point of $L$. This implies that $\widetilde{p}$ is collinear to each point of $\widetilde{L}$.
Condition~(C1) implies that there exist $p\in\widetilde{p}$, $q\in\widetilde{q}$ and $r\in \widetilde{r}$ with $p\perp q\perp r\perp p$ and so all of $p,q,r$ lie in a common member of $\mathcal{F}$, implying that $p$ is collinear to all points of the line $qr$. Hence $\widetilde{p}$ is collinear to all points of $L$.

Now assume first that $\mathcal{F}$ contains a parapolar space.  Let $\Omega_i$, $i\in I$ be such a member, let $p\in X_i$ be arbitrary and let $L\in\cL_i$ be arbitrary but such that $\delta(p,L)=2$. Suppose for a contradiction that $\widetilde{p}$ is collinear to some point $\widetilde{q}$ on $\widetilde{L}$. Let $p\perp s\perp r\in L$. Then $\widetilde{p}\perp\widetilde{q}\perp\widetilde{r}\perp\widetilde{s}\perp\widetilde{p}$, which implies by~(C1) that all of $\widetilde{p},\widetilde{q},\widetilde{r},\widetilde{s}$ contain representatives in $\Omega_i$, and hence $p$ is collinear to some point of $L$ after all, a contradiction.

Next assume all members of $\mathcal{F}$ are polar spaces. Then there exist two members $\Omega_i,\Omega_j$, $i,j\in I$, and points $p_i\in X_i$ and $p_j\in X_j$ such that $p_i$ and $p_j$ are contained in the same equivalence class $\widetilde{p}$. Choose a point $x \in X_i$ collinear to $p_i$ and a line $L$ in $\Omega_j$ not incident with $p_j$. One now easily checks that, with self-explaining notation, $\widetilde{x}$ is not collinear to any point of $\widetilde{L}$. This completes the proof of (PPS1).

\item[(PPS2)] Let $\widetilde{p}\perp\widetilde{q}\perp\widetilde{r}\perp\widetilde{s}\perp\widetilde{p}$ be a quadrangle in $\Omega$ with $\widetilde{p}$ not collinear to $\widetilde{r}$. Condition~(C1) implies that there are unique representatives $p,q,r,s$ of $\widetilde{p},\widetilde{q},\widetilde{r},\widetilde{s}$, respectively, contained in a common member $\Omega_i$ of $\mathcal{F}$, for a unique $i\in I$. Clearly, $p \perp q \perp r \perp s$ and $p$ and $r$ are not collinear. It follows that the image in $X$ of the unique symp $\xi(p,r)$ of $\Omega_i$ containing $p$ and $r$ is part of the convex closure $C$ of $\widetilde{p}$ and $\widetilde{r}$. But $C$ does not contain any further points since this would yield a circuit of length 4 and again a contradiction to Condition~(C1). 

In particular this shows that the map $p \mapsto \widetilde{p}$ is bijective when restricted to the point set of symps of $\mathcal{F}$ and hence symps occurring in members of $\mathcal{F}$ correspond bijectively with symps in $\Omega$. Hence,  the symplectic rank of $\Omega$ is the minimum of all (symplectic) ranks of the members of $\mathcal{F}$. 

\item[(PPS3)] In view of the above, this follows immediately from the fact that every line of any member of $\mathcal{F}$ is contained in a symp of that member. 
\end{itemize}
Finally, suppose that $\Omega$ has symplectic rank at least 3 (so by the above, all members of $\mathcal{F}$ have (symplectic) rank at least 3). The relation $\mathcal{R}$ is non-trivial and so there exists a class $\widetilde{p}$ with at least two elements, say $p,p'$. By Condition~(C1), $p^\perp\cap p'^\perp=\emptyset$ and moreover, no pair of points in $p^\perp\cup p'^\perp$ is in the same equivalence class. This implies that $p^\perp$ and $p'^\perp$ induce two connected components of $\Omega_{\widetilde{p}}$, so $\Omega$ is locally disconnected at~$\widetilde{p}$.  
\end{proof}

\begin{theorem}
Let $\Omega$ be a parapolar space of symplectic rank at least 3. Then, either $\Omega$ is locally connected, or $\Omega$ is isomorphic to the buttoned parapolar space $\Omega(\mathcal{F},\mathcal{R})$ where the sheet space $\mathcal{F}$ is the family of connected components of the unbuttoning $\widetilde{\Omega}$ of $\Omega$ (which are locally connected (para)polar space of symplectic rank at least 3) and (non-trivial) equivalence relation $\mathcal{R}$ defined by ``sharing the {first coordinate}''.
\end{theorem}

\begin{proof}
If $\Omega$ is not locally connected, then let $\widetilde{\Omega}$ be its unbuttoning. Let $\mathcal{F}$ be the family of connected components of $\widetilde{\Omega}$ and let $\mathcal{R}$ be the equivalence relation  on the point set of $\widetilde{\Omega}$ defined by ``sharing the {first coordinate}''. If $\mathcal{R}$ satisfies Conditions~(C1) and~(C2), it is clear that $\Omega$ is isomorphic to the buttoned geometry $\Omega(\mathcal{F,R})$ arising from $\mathcal{F}$ and $\mathcal{R}$. Now, by Lemma~\ref{universal}, $\mathcal{R}$ satisfies Condition~(C1). Moreover connectivity of $\Omega$ implies that $\mathcal{R}$ satisfies Condition~(C2). 
\end{proof}
For the applications, we need the following proposition.
\begin{prop}\label{linegraph}
Let $\Omega=(X,\cL)$ be a buttoned parapolar space (of symplectic rank at least $3$) over a  {class} $\mathcal{C}$.  Define the graph $\Gamma=(\cL,\sim)$, where $\widetilde{L}\sim \widetilde{L}'$ if $\widetilde{L}$ and $\widetilde{L}'$ are contained in a common plane. Then containment induces a bijection $\beta$ from the connected components of $\Gamma$ to the sheet space of $\Omega$.  Hence the union of all points on the lines of each connected component of $\Gamma$ is either a polar space, a locally connected parapolar space, or a self-buttoned parapolar space whose (singleton) sheet space is an element of $\cC$.
\end{prop}

\begin{proof}
Let $\widetilde{L}\in\cL$. Then, by definition, there is a line $L$ contained in some member $\Omega_*\in\cC$ contained in the sheet space such that $\widetilde{L}=\{\widetilde{p}\mid p\in L\}$. By Lemma~\ref{uniqueline}, $L$ is unique.    Now let $\widetilde{L}'$ be adjacent to $\widetilde{L}$ in $\Gamma$, and define $L'$ with respect to $\widetilde{L}'$ similarly as $L$ is defined with respect to $\widetilde{L}$. Then the first paragraph of the proof of (PPS1) in Proposition~\ref{crucial} implies that $L$ and $L'$ are both contained in $\Omega_*$. Hence $\widetilde{L}$ and $\widetilde{L}'$, and all lines in the same connected component of $\Gamma$, are contained in the image of $\Omega_*$ in $\Omega$, which is either a polar space, a locally connected parapolar space, or a buttoned parapolar space over $\{\Omega_*\}$. The proposition follows. 
\end{proof}



\subsection{Applications: extending characterizations of parapolar spaces}
\subsubsection{General discussion}
Proposition \ref{crucial} allows to remove the requirement of local connectivity in (almost) all of the characterisation results on parapolar spaces. Proposition~\ref{linegraph} allows to formulate a sufficient condition under which a characterization theorem for locally connected parapolar spaces, say satisfying the list of hypotheses (h), can be extended to locally disconnected parapolar spaces by simply choosing the sheet space in the family of all locally connected (para)polar spaces     satisfying (h).   

Indeed, we observe that almost all hypotheses of all characterization/classification results of parapolar spaces can be phrased using one and the same principle: \emph{The absence of certain configurations, that is, subgeometries}, possibly under certain additional conditions. We provide examples later; we first phrase the basic principle. 

\begin{defi}
In the present paper, a \emph{configuration} $\Sigma$ is a point-line geometry $(Y,\cM)$ with the property that every triangle of points (that is, three pairwise collinear points not on a common line) is contained in a unique projective plane. A configuration $\Sigma$ is called \emph{line-compact} if it is connected and the graph $\Gamma_\Sigma$ on the lines of $\Sigma$, with two lines adjacent if they are contained in a common plane, is connected and has diameter at most $7$. 
\end{defi}
The next proposition follows from Proposition~\ref{linegraph}, {observing that the point-diameter of a line-compact configuration is at most 3 (which follows from the restriction on the diameter of the graph $\Gamma_\Sigma$ above).} 
A \emph{full subgeometry} is a subset of points and lines such that every line of the subgeometry is a full line of $\Omega$ (however, not all pairs of points of the subgeometry which are collinear $\Omega$ are necessarily collinear in the subgeometry).


\begin{prop}\label{linecompact}
Let $\cC$ be the class of all locally connected (para)polar spaces {in which the occurrence of full subgeometries isomorphic to a member of a certain list \emph{(L)} of line-compact configurations is restricted (i.e., specific conditions are given under which they are not allowed to occur).} Then the class of locally disconnected parapolar spaces which do not admit any subspace isomorphic to a member of \emph{(L)} coincides with the class of buttoned parapolar spaces over $\cC$.    
\end{prop}     
If the hypotheses (h) of a characterization result can not be stated as the absence of line-compact configurations, then anything can happen: the locally {connected} examples could be the only examples; there could be additional locally disconnected examples, but they might be buttoned over a {differently} 
restricted set of locally connected (para)polar spaces satisfying (h), with perhaps an additional restriction on the buttoning, for instance on the graph in Condition (C2). We will illustrate this with examples. 
In any case, these examples will show that, in general, it would still be worthwhile to spend a few words on the locally disconnected case every time a new characterization/classification result is discovered, since it might not be so straightforward to apply Shult's principles, alluded to in the introduction of this appendix. These principles state that the condition of local connectivity can be deleted provided  that either 

\begin{compactenum}
\item[(P1)] one adds in the conclusion buttoned parapolar spaces over the class of locally connected (para)polar spaces satisfying all other requirements, that is, the class of (para)polar spaces satisfying the given conditons (other than local connectivity) is closed under buttoning and unbuttoning; 

\item[(P2)] If (P1) does not hold, then one can extend the classification in question to the buttoned parapolar spaces by requiring that the given conditions must hold in every member of the sheet space.   
\end{compactenum}
We will refer to the principles (P1) and (P2) below. It is clear that (P1) requires some checking, but (P2) is rather ``an easy way out'' if (P1) is not (obviously) applicable. The goal of this revision is to show that, if (P1) fails, then there are sometimes more subtle and more efficient ways to generalize the given characterization than suggested by (P2), see for instance our comment concerning the $(-1)$-lacunary parapolar spaces. 

\subsubsection{Lacunary parapolar spaces}
In \cite{MinusOnePaper} and the present paper, all locally connected lacunary parapolar spaces of  symplectic rank at least  $3$ are classified,  and also all strong lacunary parapolar spaces of arbitrary symplectic rank. 
Note that by Lemma~\ref{index0strong}  all $0$-lacunary parapolar spaces are strong and hence, if the symplectic rank is at least $3$, they are locally connected. Also, the condition of being $k$-lacunary, $k\geq 1$ is equivalent to the absence of the line-compact configurations consisting of two polar spaces of rank $\geq k+1$ intersecting in exactly a singular $k$-space, subject to the condition that these polar spaces are symps.   

This implies the following reduction, which conforms to the principle (P1).

\begin{theorem}\label{laclac}
Let $k\geq 1$. Then a locally disconnected parapolar space $\Omega$ of symplectic rank at least $3$ is $k$-lacunary if and only if it is a buttoned parapolar space over the family $\cC_k$ of all polar spaces of rank at least $\max\{3,k+1\}$ and locally connected $k$-lacunary parapolar spaces of symplectic rank at least $3$.
\end{theorem}

%

Proposition~\ref{linecompact} is not applicable to the $(-1)$-lacunary condition. Indeed, in terms of nonexistence of a certain configuration, this is just the disjoint union of two polar spaces, which must be symps, and this is not line-compact (even not connected). However, the locally disconnected case has already been settled in \cite{MinusOnePaper}, and it shows indeed that it is not just a case of restricting the sheet space. Indeed,  Theorem~3.1 of \cite{MinusOnePaper} implies that, as soon as a line is contained in at least two symps, then $\Omega$ is locally connected. Also, Theorem~3.2 of \cite{MinusOnePaper} can be restated as \begin{quote} \em a locally disconnected parapolar space of symplectic rank at least $3$ is $(-1)$-lacunary if and only if it is a buttoned parapolar space over the class of all polar spaces of rank at least $3$, satisfying additionally that the graph in \emph{Condition (C2)} is complete.\end{quote}In particular, this classification does not fit in any of Shult's principles (P1) or (P2).  

\begin{remark} Note that (C1) implies that an equivalence class of $\mathcal{R}$ cannot contain two points of the same member of $\mathcal{F}$ if that member's diameter is smaller than 5. This situation in particular applies if all members of $\mathcal{F}$ are ordinary polar spaces; which is always the case when $k=-1$ and when $k \geq 6$; indeed, if $k\geq 6$, then there are no locally connected $k$-lacunary parapolar spaces due to the main results of the present paper. The latter is an example that locally connected parapolar spaces satisfying certain hypotheses might not exist although locally disconnected ones do exist. This fits into Principle (P1).
\end{remark}

Here is an application of the buttoned geometries over a family containing polar spaces of rank 2, outside the realm of the local connectivity theory {which required rank at least 3}.  
  
It is shown in Lemma~\ref{1,minrank>2} that strong $1$-lacunary parapolar spaces with minimum symplectic rank 2 do not exist. We can however drop the assumption of strongness and classify all $1$-lacunary parapolar spaces with minimum symplectic rank 2. 

\begin{theorem}\label{1lacnonstrong}
A parapolar space $\Omega$ of minimum symplectic rank $2$ is $1$-lacunary if and only if it is a buttoned parapolar space over the family $\cC_1'$ of all polar spaces of rank at least $2$ and locally connected $1$-lacunary parapolar spaces of symplectic rank at least 3.
\end{theorem}

\begin{proof}
Let $\Omega=(X,\cL)$ be a $1$-lacunary parapolar space with minimum symplectic rank $2$.  Set \[X'=\{x\in X: x\mbox{ is contained in a symp of rank at least }3\},\] \[\cL'=\{L\in\cL:L\mbox{ is contained in a symp of rank }3\}.\] 
Then each connected component $\omega$ of $\Omega'=(X',\cL')$ is a (para)polar space with symplectic rank at least $3$ (observe that a line of $\mathcal{L}'$ cannot be contained in a symp of rank 2 as well, since two symps cannot share exactly a line by assumption). Moreover, if $\omega$ is a parapolar space, it is $1$-lacunary. Considering the disjoint union of all connected components of $\Omega'$ and all symps of $\Omega$ of rank 2 and then identifying all points which were the same in $\Omega$ (via an obvious equivalence relation $\mathcal{R}$), we see that we obtain  Construction~\ref{buttoned} (note also that no connected component is self-buttoned as Proposition~\ref{prop12} implies that $1$-lacunary parapolar spaces of symplectic rank at least 3 have diameter at most $3$). 
{Hence, $\Omega$ is indeed a buttoned parapolar space over the family $\cC_1'$.}
The converse is obvious.
\end{proof}

We now quickly review some other results in the literature and formulate the corresponding theorems without the assumption of being locally connected. 

\subsubsection{On a theorem of Shult}
In Theorem 19 of \cite{Shu:05}, Shult classifies all locally connected parapolar spaces of symplectic rank at least $4$ possessing a collection $\cM$ of maximal singular subspaces satisfying these properties:
\begin{compactenum}[$(i)$]
\item Every plane lies in a member of $\cM$.
\item If $M \in \cM$, and $x$ is a point not in $M$ then $x^\perp\cap M$ is either empty, a single
point, or a $3$-dimensional projective space.
\end{compactenum}
The conclusion of Shult's theorem are the parapolar spaces $\mathsf{E}_{n,1}(*)^h$, $n\geq 6$.

The first condition is trivial: the configuration is just a plane, and the additional condition is that it is not contained in a member of $\cM$ (i.e., there cannot be planes not occurring in a member of $\cM$). The second condition is also easy to rephrase as {the absence of the following line-compact configurations}: 
{the union of two distinct maximal singular subspaces $U$ and $M$, intersecting each other in a non-empty subspace of dimension not equal to $0$ or $3$, with the additional condition that $M$ belongs to $\mathcal{M}$.} 
Hence similarly as Theorem~\ref{laclac}, and combining it with Theorem~1 of \cite{Shu:05}, we have the following complement to Shult's result, conforming to Principle (P1).
\begin{theorem}\label{Shultagain}
Let $k\geq 3$. Let $\Omega$ be a locally disconnected parapolar space of symplectic rank at least $4$ possessing a collection $\cM$ of maximal singular subspaces satisfying these properties:
\begin{compactenum}[$(i)$]
\item Every plane lies in a member of $\cM$.
\item If $M \in \cM$, and $x$ is a point not in $M$, then $x^\perp\cap M$ is either empty, a single
point, or a $k$-dimensional projective space.
\end{compactenum}
Then $\Omega$ is a buttoned parapolar space over the family $\cC_k^*$ consisting of polar spaces of rank $k+2$, if $k>3$, and polar spaces of rank $5$ and the locally connected parapolar spaces $\mathsf{E}_{n,1}(*)$, $n\geq 6$, if $k=3$. 
\end{theorem}

To show that care is needed, we mention the following special case of Shult's theorem, and note that the straightforward similar special case of Theorem~\ref{Shultagain} is \emph{not} true.

\begin{corollary}
 Let $\Omega$ be a locally connected parapolar space of symplectic rank at least $4$ possessing a collection $\cM$ of maximal singular subspaces satisfying these properties:
\begin{compactenum}[$(i)$]
\item Every plane lies in a member of $\cM$.
\item If $M \in \cM$, and $x$ is a point not in $M$ then $x^\perp\cap M$ is either a single
point, or a $3$-dimensional projective space.
\end{compactenum}
Then $\Omega$ is the strong parapolar space $\mathsf{E_{6,1}}(\K)$, for some field $\K$.   
\end{corollary}

Although the conditions of this corollary assume local connectivity, and the proof requires it, it is clear that no buttoned parapolar space can ever satisfy the two conditions $(i)$ and $(ii)$ simultaneously. Indeed, if $\Omega$ is locally disconnected at some point $p$, then in one component of $\Omega_p$, we find a member of $\cM$, but in the other component, we find a point not collinear to $p$, a contradiction. So, after all, the local connectivity assumption can be deleted without any consequence on the spaces in the conclusion. Note that the new condition $(ii)$ cannot be stated as the absence of a line-compact configuration, since we now have to add the configuration consisting of a member of $\cM$ and a disjoint point, which is not even connected. In general, it is tempting to state that, if a ``reasonable'' configuration appearing in the ``absence''-conditions is disconnected, then the polar spaces in the conclusions must automatically be locally connected. However, this is not true for $(-1)$-lacunarity, although we could argue that the corresponding configuration is not ``reasonable'' since it comprises two symps, and these can each fill a whole component (sheet).  

This example shows that, although formally obeying the principles (P1) and (P2) since no locally disconnected ones satisfy the hypotheses, one cannot blindly apply them. 

\subsubsection{On a theorem of Cohen-Cooperstein}

As a final example we mention Shult's updated version in \cite{Shu:12} of a theorem of Cohen and Cooperstein \cite{Coco}, restricted to symplectic rank 3. One considers locally connected parapolar spaces of uniform symplectic rank $3$ such that the intersection $x^\perp\cap\xi$, with $\xi$ a symp not containing the point $x$, is never a line. This translates easily in the absence of a line-compact configuration and hence, completely similarly as above, the locally disconnected case is just the buttoning over the family consisting of all locally connected parapolar spaces satisfying this requirement, together with all polar spaces of rank $3$.

{\bf Acknowledgement} All four authors are very grateful to The University of Auckland Foundation which awarded the third author a Hood Fellowship Grant nr 3714543, during his stay as Hood Fellow when also authors 1 and 4 visited the University of Auckland. The majority of this paper was written in this period. The fourth author would also like to thank the fund Professor Frans Wuytack for their support.

\end{document}